\documentclass{amsart}
\usepackage{txfonts}
\usepackage{amscd}
\usepackage{cite}
\usepackage{epsfig}
\usepackage{verbatim}
\usepackage{mathdots}
\usepackage{amssymb}
\usepackage{amsfonts}
\usepackage{amsbsy}
\usepackage{graphicx}
\usepackage[usenames]{color}
\usepackage{cancel}
\usepackage{multirow}
\usepackage{empheq}
\usepackage{cases}
\usepackage{enumitem}

\newtheorem{theorem}{Theorem}[section]
\newtheorem{corollary}[theorem]{Corollary}
\newtheorem{proposition}[theorem]{Proposition}
\newtheorem{lemma}[theorem]{Lemma}
\newtheorem{question}[theorem]{Question}
\theoremstyle{definition}
\newtheorem{definition}[theorem]{Definition}
\theoremstyle{remark}
\newtheorem{remark}[theorem]{Remark}
\newtheorem*{remark*}{Remark}
\newtheorem{example}[theorem]{Example}
\numberwithin{equation}{section}

\newcommand{\R}{{\mathbb R}}

\newcommand{\dist}{\operatorname{dist}}

\newif\ifdraft
\draftfalse
\drafttrue
\ifdraft
\else
\fi
\begin{document}	
	\title{Arithmetic product theorems on self-similar sets of uniqueness and homogeneous Moran sets}
		
	\author[Y.-S. Fu]{Yan-Song Fu}
	\address{Department of  Mathematics   \\ China University of Mining and Technology (Beijing) \\ Beijing, 100083\\P. R. China.\\ https://orcid.org/0000-0003-4687-5729}
	\email{yansong\_fu@126.com }
		
	\author[X.-G. Wang]{Xu-Guang Wang}
	\address{Department of  Mathematics   \\ China University of Mining and Technology (Beijing) \\ Beijing, 100083\\P. R. China.\\ https://orcid.org/0009-0002-3433-0713}
	\email{wxg1211@gmail.com}
	
	\author[M. Zhu]{Meng Zhu} \address{School of Mathematics and Statistics   \\ Beijing Jiaotong University \\ Beijing, 100044 \\P. R. China.\\ https://orcid.org/0000-0001-7000-8823} \email{zhumengsx86@163.com;}
\begin{abstract}
Varj\'u and Yu have given a complete description of non-singleton self-similar sets of uniqueness on the real line $\R$. In this paper, we first prove an arithmetic product theorem for self-similar $U$-sets $E:$ if \(F\subset\mathbb R\) is compact, then $F\cdot E$ is an $M_0$-set if and only if $F$ is uncountable. We next show that this conclusion remains valid for bounded integer-based homogeneous Moran sets. As an application, for every $s\in(0,1]$, there exists a bounded integer-based homogeneous Moran set of uniqueness  $K$   such that $\dim_H K=s$ and the self-product $K\cdot K$ is an $M_0$-set. In sharp contrast,  for every $s\in[0,1/2]$ we can construct an  unbounded integer-based homogeneous Moran $U$-set $E$ with $\dim_H E=s$ such that both $E$ and $E\cdot E$ are $U$-sets, and $1/2$ is the optimal threshold for $U$-sets to be closed under the arithmetic  self-product.  Moreover, we show that if $E\subset[0,1]$ is a compact set  whose self-product $E\cdot E$ is contained in an $H$-set, then $\underline\dim_B E\le1/2$.
	\end{abstract}

	\def\subjclassname{2020 Mathematics Subject Classification}
	\subjclass[2020]{Primary 42A63, 28A80, 42A38, 11J72.}
	
	
	\maketitle
	
	\section{Introduction}
	In classical Fourier analysis,  a set $E\subset [0,1]$ is called a set of \textit{uniqueness}, or $U$-set, if every trigonometric series converging to $0$ outside $E$ vanishes identically, i.e., $\sum_{n\in\mathbb Z}c_ne^{2\pi i nx}=0$ for all $x\in[0,1]\setminus E$ implies that $c_n=0$ for all $n$. A set $E\subset \R$ is a $U$-set if its image in $\R/\mathbb{Z}=[0,1]$ is. The rest of sets are called the sets of \textit{multiplicity}, or $M$-sets. In particular, a set $E\subset \R$ is called an $M_0$-set, or $M$-set in the \textit{narrow sense}, if $E$ supports a Borel probability measure $\mu$ such that the Fourier transform
	\begin{equation*}
		\widehat{\mu}(u)=\int e^{-2\pi i ut}d\mu(t)
	\end{equation*}
	satisfies that $ \widehat{\mu}(u)\to0$ as $|u|\to\infty.$ Such measures are called \textit{Rajchman measures}.  Therefore, constructing Rajchman measures provides a way to verify that $E$ is an  $M_0$-set. Of course, there are examples of $M$-sets which are not $M_0$-sets; see  Piatetski-Shapiro \cite{Piatetski-Shapiro1952,Piatetski-Shapiro1954}. If a set $E$ is not an $M_0$-set, it is called a $U$-set in the \textit{broad sense}, or $U_0$-set.

	The research of $U$-sets and $M$-sets dates back to the work of Riemann and Cantor in the mid-19th century on the uniqueness of the representation by trigonometric series.  By  Riemann's theory on trigonometric series (see, e.g., \cite[Chapter IX]{Zyg1959}), every countable set is a $U$-set and every set with positive Lebesgue measure is an $M$-set.  However, the problem of the classification of sets of Lebesgue measure zero into $U$-sets and $M$-sets is very difficult and far from solved. We can refer to the books of Zygmund \cite{Zyg1959}, Salem \cite{Salem1963}, Kahane and Salem \cite{KahaneSalem1963}, and Kechris and Louveau \cite{KechrisLouveau} for the classical theory about $U$-sets and $M$-sets.
	
The question we investigate in this paper is related to the research  on the Lebesgue measure of arithmetic product of singular fractal sets and Fourier decay of multiplicative convolutions.  By $A\cdot B$ we denote the {\it arithmetic product} of subsets \(A,B\subset\mathbb R\), i.e., 	\(A\cdot B=\{ab:a\in A,\ b\in B\}\).   Athreya, Reznick and Tyson \cite{ART2019} proved that, for the middle-third Cantor set $C_{1/3}$,  the arithmetic product set $C_{1/3}\cdot C_{1/3}$ has non-trivial intervals and its Lebesgue measure is between  $17/21$ and $8/9$, yielding that $C_{1/3}\cdot C_{1/3}$ is an $M_0$-set. Other related arithmetic product sets $E\cdot F$ having positive Lebesgue measure have also been obtained under the assumption that the sets  $E,F$ with `large' Hausdorff dimension, see \cite[Theorem 1.9]{VY2026} and \cite[Theorem 1.1]{GJXZ2021}. In particular, if the arithmetic product set $E\cdot F$ has nonempty interior, then it is an $M_0$-set. The arithmetic product set carries a natural measure, known as multiplicative convolutions, see Bourgain \cite{Bourgain2010} and  the references cited therein. Recently,  Orponen, Saxc\'e and Shmerkin used the \emph{\(s\)-energy}, for \(s \in (0,1]\), of finite Borel measure  \(\mu\) on \(\mathbb{R}\),  defined as
	\begin{equation*}
		I_s(\mu) := \iint |x - y|^{-s} \, d\mu(x) \, d\mu(y),
	\end{equation*}
	to prove the following result:
	\begin{theorem}[\textbf{Fourier decay of multiplicative convolutions}]\label{thm.OSS}
		Let \(n\ge 2\), and \(\{s_j\}_{j=1}^n \subset (0,1]\) such that \(\sum s_j > 1\). There exists \(\tau = \tau(n,\{s_j\}) > 0\) such that the following holds. Let \(\mu_1,\dots,\mu_n\) be Borel probability measures supported on \(E_1, \dots,E_n\subset \mathbb R\) respectively, such that \(I_{s_j}(\mu_j) < +\infty\). Then there is \(C = C(\{\mu_j\}) > 0\) such that
		\begin{equation*}
			\left| \int e^{-2\pi i \xi x_1\cdots x_n} \, d\mu_1(x_1)\dots d\mu_n(x_n) \right| \le C \cdot |\xi|^{-\tau}, \quad \xi \in \mathbb R.
		\end{equation*}Consequently,  the arithmetic product $E_1\cdot E_2\cdots E_n$ is an $M_0$-set.
	\end{theorem}
	
	Therefore, the preceding results naturally give rise to the following question:
	\begin{question}\label{ques}
		Given two sets $E,F\subset \R$ with Lebesgue measure zero, how does one determine whether the arithmetic product set $E\cdot F$ is an $M_0$-set  or a $U$-set?
	\end{question}
	
	In this paper, we will address this question by fixing one factor to be self-similar sets or homogeneous Moran sets, which are important subjects of study in fractal geometry, geometric measure theory, harmonic analysis, dynamical system and related fields.

	\subsection{The self-similar sets of uniqueness and multiplicity}\label{sec.1}
	
Recall that, for an \textit{iterated function system} (IFS) written as
 	\begin{equation}\label{eq.IFS}
\Phi:=\{f_i(x)=r_i x+a_i:1\le i\le k\},\quad k\ge2,\quad 0<r_i<1, \quad a_i\in\R,
	\end{equation}
	Hutchinson \cite{Hut1981} showed that there is a unique nonempty compact set \(E\), called a \textit{self-similar set}, such that
	\[E=\bigcup_{i=1}^{k}f_i(E);\]
	moreover, for a probability vector \((p_1,\ldots,p_k)\) with strictly positive entries such that  \(\sum_{i=1}^k p_i=1\), there is a unique probability measure \( \mu\), called a \emph{self-similar measure}, supported on the set $E$, such that
	\begin{equation}\label{eq.meas}
		\mu = \sum_{i=1}^k p_i\mu\circ f_i^{-1}.
	\end{equation}
Here and below, \(\mu\circ f_i^{-1}\) denotes the pushforward of \( \mu \) under \(f_i\), i.e.,  \((\mu\circ f_i^{-1})(E) = \mu(f_i^{-1}(E)) \) for all Borel sets \(E\). The self-similar IFS is \textit{homogeneous} if all $r_i$ are equal to each other, i.e., $r_i=r$ for all $i.$
	
	The most well-known self-similar set is the  \emph{middle-$r$ Cantor set} \(C_r\), where \(0<r< \frac{1}{2}\).  Remember  that \(C_r\) is obtained by repeatedly removing a proportion \(0<1-2r <1\) from the middle of intervals, starting with the unit interval $[0,1]$, whose Hausdorff and box dimensions are \(-\log 2/\log r\). The classical middle-third Cantor set is obtained by taking $r=1/3$. The celebrated work of Salem and Zygmund \cite{SalemZygmund1955} (also see \cite[Chapter VI]{Salem1963}) showed that \(C_r\) is a set of uniqueness if and only if $r^{-1}$ is a Pisot number, obtaining the connection between harmonic analysis and number theory.  A positive number $\theta>1$ is a \textit{Pisot number} if it is an algebraic integer whose Galois conjugates have modulus strictly less than $1$.
	
	The breakthrough result of Li and Sahlsten \cite{LiSahlsten2022}  showed that if the self-similar set $E$ generated by the IFS in \eqref{eq.IFS} is non-singleton and \(\log r_j / \log r_\ell\) is irrational for some \( j \neq \ell\), then the self-similar set $E$ is an $M_0$-set, by showing that every  self-similar measure $\mu$ in \eqref{eq.meas} supported on $E$ is a Rajchman measure. This result reduces the problem of characterizing when a self-similar set is a \(U\)-set or an \(M\)-set, and when a self-similar measure is non-Rajchman, to the case \(r_j=r^{l_j}\) for some \(r\in (0, 1)\) and \( l_j \in \mathbb{Z}_{>0}\) with \(\gcd(l_1,\ldots,l_k)=1\). In this case, Br\'{e}mont \cite{Bremont2021} showed that if \( \mu\) is a non-singleton self-similar measure that is not Rajchman, then \( r^{-1}\) is a Pisot number and the IFS can be conjugated by a suitable similarity to a form such that \( a_i \in \mathbb{Q}(r)\) for all \(i\). Finally, Varj\'u and Yu \cite{VarjuYu2022}  completely identify all self-similar sets which are sets of uniqueness, which is stated as follows.
	\begin{theorem}[Varj\'u and Yu {\cite[Theorem 1.4]{VarjuYu2022}}]
		\label{thm:intro-varju-yu-form}
		Let \(E\subset\mathbb R\) be a non-singleton self-similar set.  Then \(E\)
		is a set of uniqueness if and only if it has zero Lebesgue measure and its
		underlying IFS can be conjugated by a suitable similarity to the form
		\[f_i(x)=r^{\ell_i}x+a_i,\qquad i=1,\dots,k,	\]
		where \(0<r<1\), \(\ell_i\in\mathbb Z_{>0}\),
\(\gcd(\ell_1,\dots,\ell_k)=1\), \(\lambda=r^{-1}\) is a Pisot 		number, and \(a_i\in\mathbb Q(\lambda)\) for all \(i\).
	\end{theorem}
	Many recent papers also have provided characterizations of the Fourier decay of higher-dimensional  self-similar measures, self-affine measures and other measures, thereby contributing to the characterization of $M$-sets. For our purpose, we do not aim to give a comprehensive survey, the interested readers can refer to  \cite{LS2020,AHW2022,R2022,S2022} and references therein.

	\subsection{The homogeneous Moran sets  of uniqueness and multiplicity}	
	Let \(\{m_k\}_{k\ge1}\) be a sequence of positive integers such that $m_k\ge2$, and let $D_k\subset\{0,1,\ldots,m_k-1\}$ be a finite subset such that $2\le\#D_k\le m_k$ for each \(k\ge 1\). Let $Q_k=m_1\cdots m_k$.  Then the associated \textit{integer-based homogeneous  Moran set} is the unique compact set $K$ such that
	\begin{equation}\label{1.12}
		K=\left\{\sum_{k=1}^{\infty} \frac{d_k}{Q_k}:d_k\in D_k\right\}=\bigcap_{n=1}^{\infty}\bigcup_{(d_1,\ldots,d_n)\in D_1\times\cdots\times D_n}I({d_1,\cdots, d_n}),
	\end{equation}
	where the $n$ order fundamental interval is defined by
	\[I({d_1,\cdots, d_n}):=\left[\sum_{k=1}^n\frac{d_k}{Q_k},\ \sum_{k=1}^n\frac{d_k}{Q_k}+\frac1{Q_n}\right].\]
It is clear that $K$ is an uncountable perfect set.	For simplicity, we call the homogeneous  Moran set $K$ the \textit{bounded integer-based} if $\sup_k m_k<\infty$, and \textit{unbounded integer-based} otherwise.
	The homogeneous  Moran set $K$ naturally supports a Borel probability measure which is the infinite convolution
	\begin{equation}\label{eq.4.1}
\mu:=\overset{\infty}{\underset{k=1}{\scalebox{2}{\(\ast\)}}} \left(\frac{1}{\#D_k}\sum_{d_k\in D_k} \delta_{\frac{d_k}{Q_k}}\right).
	\end{equation}
Here and below, $\delta_a$ denotes the Dirac measure concentrated at the point $a$, the asterisk symbol $\ast$ denotes the convolution of two measures, and $\#D$ denotes the cardinality of a finite set $D$. The Fourier transform of the measure $\mu$ is the infinite product
	\[\widehat\mu(u)=\int e^{-2\pi iu\cdot x}\,d\mu(x)=\prod_{k=1}^{\infty}P_{D_k}\left(\frac{u}{Q_k}\right).\]
	Here and below, for a digit set $D$, the symbol $P_D$ denotes the mask polynomial \[P_D(u)=\frac1{\#D}\sum_{d\in D}e^{-2\pi i d u}.\]
	
The study on Homogeneous Moran sets and measures has a long history in harmonic analysis, number theory and fractal geometry; see  \cite{KahaneSalem1963,FengWenWu1997,Kaufman1995}. For the early classification of homogeneous Moran sets into
$U$-sets and $M$-sets, see \cite[Chapters 5 and 6]{KahaneSalem1963}. 
Lai \cite{Lai2017} proved that if $\sup_{n\geq 1} \frac{\max D_n}{m_n} \leqslant L < 1/2,$ then \(K\) defined in \eqref{1.12} will be a $U$-set. Very recently, Lai and Xie \cite[Proposition 6.1]{LaiXie2026} showed that there is an {\it unbounded integer-based} homogeneous Moran set of uniqueness has Hausdorff dimension $1$, while the self-similar set of uniqueness must have Hausdorff dimension strictly less than $1$ (\!\!\cite[Proposition 1.5]{LaiXie2026}). In fact, the conclusion of \cite[Proposition 6.1]{LaiXie2026} holds for all $s\in(0,1]$ by modifying its proof; see \cite[Subsection 6.1]{LaiXie2026}.  Furthermore,   Lai, Shi and Xie \cite[Theorem 1.1]{LaiXie20261} explicitly proved a  {\it dimension intermediate value theorem}  for $M_0$-set, that is,  for every $s\in (0,1]$ there is a homogeneous Moran set $E$   of multiplicity of Hausdorff dimension $s.$

	\subsection{Main results}
	
	Our first result answers Question \ref{ques} concerning $M_0$-set affirmatively for every non-singleton self-similar \(U\)-set on $\R$ described in Theorem \ref{thm:intro-varju-yu-form}.
	\begin{theorem}\label{thm:intro-self-similar-product}
		Let \(E\subset\mathbb R\) be a non-singleton self-similar set of uniqueness, and let \(F\subset\mathbb R\) be a compact set. Then \[F\cdot E \text{ is an } M_0\text{-set } \quad\text{if and only if}\quad	F\text{ is uncountable}.\]
	\end{theorem}
	To the best of our knowledge, Kaufman \cite{Kaufman1993} first proved that Theorem \ref{thm:intro-self-similar-product} holds for the middle-third Cantor set $C_{1/3}$, which was used to show that the class $M_0$  of closed sets of multiplicity is analytic but not Borel. In this sense,  Theorem \ref{thm:intro-self-similar-product} extends   Kaufman's result to all  non-singleton self-similar set of uniqueness on $\R$. Moreover, we give Example \ref{ex:pisot-interval} to show that the uniqueness assumption in  Theorem \ref{thm:intro-self-similar-product} cannot be dropped.

	\noindent\textbf{The sketch of proof of Theorem \ref{thm:intro-self-similar-product}.}  The necessity is a direct application of two celebrated theorems. In fact, it is well-known that Marcinkiewicz-Zygmund theorem (see, e.g., \cite[Theorem 1]{MZ1937}, \cite[p. 71, Theorem 2]{KechrisLouveau}) showed that if \(a\ne0\), then \(aE\) is a \(U\)-set if and only if \(E\) is a \(U\)-set. Bari  (\!\cite[p. 10]{Bari1949}, also see \cite[p. 41, Theorem 1]{KechrisLouveau}) proved that the countably many  union of closed \(U\)-sets is also a \(U\)-set. Therefore, if $F$ is a countable set, then  the set \(F\cdot E=\bigcup_{a\in F}aE\) is a  \(U\)-set, which yields that the necessity of Theorem \ref{thm:intro-self-similar-product} follows. The sufficiency is established provided there exists a probability measure $\mu$ supported on $E$  satisfying  the following two-scale Fourier estimate
	\begin{equation*}
		\widehat\mu(u)\,\widehat\mu(\tau u)\to0
		\qquad(|u|\to\infty)
	\end{equation*}
	for every \(\tau\) outside some countable set (Lemma \ref{lem:product-kaufman}).  Once the above two-scale Fourier estimate holds for $\mu$, pushing the product measure on $F\times E:=\{(x,y): x\in F, y\in E\}$ forward under the mapping  \((s,x)\mapsto sx\) yields a Rajchman measure supported on \(F\cdot E\), which shows that \(F\cdot E\) is an $M_0$-set. Thus the key input is the two-scale Fourier estimate together with the fact that its exceptional set is countable. This allows the product argument to apply to every uncountable compact set $F$, including those with Hausdorff dimension zero, i.e., $\dim_HF=0$.

	Since each non-singleton self-similar set is uncountable and $M_0$-set is closed under arithmetic products with arbitrary nonzero sets, Theorem \ref{thm:intro-self-similar-product} immediately yields   the following conclusion.	
	\begin{corollary}\label{cor:intro-finite-self-similar-products}
			Let \(n\ge2\), and let \(E_1,\ldots,E_n\subset\mathbb R\) be non-singleton self-similar sets.  If  at least one of   \(E_1,\ldots,E_n\) is a $U$-set, then \(E_1\cdot E_2\cdots  E_n\) is an $M_0$-set.
	\end{corollary}

We now study  Question \ref{ques} for the  homogeneous Moran sets. We  answer  Question \ref{ques} for the  bounded integer-based case  by establishing an analogue of Theorem \ref{thm:intro-self-similar-product}.

	\begin{theorem}\label{thm:intro-bounded-integer-moran}
		Let \(K\) be the homogeneous Moran set as in \eqref{1.12} such that there exists a constant $B$ such that \(2\leq\#D_k\le m_k\le B<\infty\) for every $k$. If \(F\subset\mathbb R\) is an uncountable compact set, then \(F\cdot K\) is an \(M_0\)-set. If, in addition, \(K\) is a set of uniqueness, then for every compact set \(F\subset\mathbb R\), \[F\cdot K\text{ is an $M_0$-set } \quad\text{if and only if}\quad F\text{ is uncountable}.\]
	\end{theorem}
	
As an application of Theorem \ref{thm:intro-bounded-integer-moran}, we establish a {\it dimension  intermediate value theorem} for $U$-sets whose self-product is an $M_0$-set, which strengthens the related results of Lai and Xie  \cite[Proposition 6.1]{LaiXie2026} mentioned before. 
	\begin{theorem}\label{rem:bounded-all-dimensions}
		For every $s\in(0,1]$, there exists a bounded integer-based homogeneous Moran set $K$ such that $K$ is a $U$-set, and  \[\dim_H K=s \quad\text{and}\quad K\cdot K\text{ is an }M_0\text{-set}.\]
\end{theorem}

 However, the phenomena described in Theorems \ref{thm:intro-self-similar-product} and \ref{thm:intro-bounded-integer-moran} may not hold in the case of unbounded integer-based homogeneous Moran sets.  We next answer Question \ref{ques} for a particular homogeneous Moran measures.

	\begin{theorem}\label{thm:unbounded-m0-product}
		Let \(K\) be the homogeneous Moran set as in \eqref{1.12} such that $D_n=\{0,1,\ldots,N_n-1\}$, where $2\le N_n\le m_n$ and $m_n\to\infty$.  If  $\inf_n \frac{N_n}{m_n}=0$, then the following two statements holds.
\begin{enumerate}[label=\textup{(\alph*)}]
  \item There exists a Cantor set $F\subset[\frac14,\frac12]$ with $\dim_HF=0$ such that $F\cdot K $ is  a $U$-set.
  \item For every $s\in(\frac12,1]$, there exists an unbounded integer-based homogeneous Moran set $K$  with $\dim_H K=s$ such that $K$ is a $U$-set, and $K\cdot K$ is an $M_0$-set
\end{enumerate}
 Furthermore,  the following two statements are equivalent:
		\begin{enumerate}[label=\textup{(\roman*)}]
			\item \label{item:positive-relative-density}
		The relative digit sizes have a positive lower bound:
			\begin{equation*}
				\inf_{n\ge1}\frac{N_n}{m_n}>0.
			\end{equation*}
			\item \label{item:universal-product}
For every uncountable compact set $F\subset\mathbb R$,	the set $F\cdot K $ is an $M_0$-set.
		\end{enumerate}
		If (i)  holds and  there exists a constant  $0<L<1$ such that $N_n\le Lm_n$ for every $n$, then $K $ is an $U$-set with $\dim_HK =1$ and, for every compact set $F\subset\mathbb R$,
		\begin{equation}\label{eq.equiv}
			F\cdot K \text{ is an }M_0\text{-set}		\quad\text{if and only if}\quad F\text{ is uncountable}.
		\end{equation}
	\end{theorem}
	
	Finally, in terms of Hausdorff and lower box dimensions, we provide an answer to Question \ref{ques} concerning $U$-sets.
	Recall that the {\it lower box dimension}  of a bounded set	$E\subset\mathbb R$ is defined by \[\underline\dim_B E=\liminf_{\delta\to0} \frac{\log N(E,\delta)}{\log(1/\delta)},\] 	where $N(E,\delta)$ is the least number of sets of diameter at most $\delta$ needed to cover $E$. A compact set  \(E\subset[0,1]\) is an $H$-set if  there exist a strictly increasing sequence of positive integers \(n_j\to\infty\) and a nonempty open interval \(I\subset(0,1)\) such that $\{n_jx\}\notin I  (x\in E,\ j\ge1).$
	
	\begin{theorem}\label{thm.1.8}
		\label{thm:sharpness}
		\begin{enumerate}[label=\textup{(\alph*)}]
			\item For every $s\in[0,1/2]$, there exists an unbounded integer-based homogeneous Moran set $E\subset[0,1]$  such that $E$ is a $U$-set,  $E\cdot E$ is an $H$-set (hence a $U$-set), and \[\dim_HE=s,\qquad \dim_H(E\cdot E)\ge s.\]
			\item Conversely, if $E\subset[0,1]$ is compact and $E\cdot E$ is contained in an $H$-set, then
			\[\dim_HE\ \le\ \underline\dim_BE\ \le\ \tfrac12.\]
			More generally, if $E_1,E_2\subset[0,1]$ are nonempty compact sets and $E_1\cdot E_2$ is contained in an $H$-set, then 	$\underline\dim_BE_1+\underline\dim_BE_2\le1$; the displayed bound follows by taking $E_1=E_2=E$.
		\end{enumerate}
	\end{theorem}
 According to Frostman's lemma \cite[Theorem 2.8]{Mat2015}, the Hausdorff dimension of a Borel set $E \subset \mathbb{R}^d$ can be expressed as $$\dim_{H} E = \sup \{s\leq d: \text{there is } \mu \in \mathcal{M}(E) \text{ such that } I_s(\mu) < \infty\},$$ where   $\mathcal{M}(E)$ is the family of Radon measures $\mu$ with compact support contained in  $E$ and $0<\mu(\mathbb{R}^d)<\infty.$
	Thus, by Theorem \ref{thm.OSS}, if the compact sets 	$E_1,E_2\subset\mathbb R$ satisfy $\dim_HE_1+\dim_HE_2>1$, then  $E_1\cdot E_2$ is an $M_0$-set. In other words,
	\begin{equation*}
		E_1\cdot E_2\text{ is not an }M_0\text{-set}		\quad\Longrightarrow\quad \dim_HE_1+\dim_HE_2\le1.
	\end{equation*}
	In particular, $E\cdot E$ being a $U$-set implies $\dim_HE\le1/2$. With this in mind,  we can make the following remarks regarding  our arithmetic product theorems on homogeneous Moran sets.
	\begin{itemize}
\item  As a consequence of Theorem \ref{thm.OSS} and Theorem \ref{thm.1.8}(a), we have
		\begin{equation*}
			\max\bigl\{\dim_HE:E\subset[0,1]\text{ is compact and both }E\text{ and } E\cdot E\text{ are }U\text{-sets}\bigr\}=\frac12,
		\end{equation*}
\begin{equation*}
			\min\bigl\{\dim_HE:E\subset[0,1]\text{ is compact and both }E\text{ and } E\cdot E\text{ are }U\text{-sets}\bigr\}=0.
		\end{equation*}
This, along with Theorem \ref{thm:unbounded-m0-product}(b), implies that the value $\frac12$ is the optimal threshold for $U$-sets to be closed under arithmetic self-product.    Moreover, if $\dim_H E\in (0,1/2)$, then for any $s<\dim_H E$  we have $2s<1$. In this case, the hypothesis of Theorem \ref{thm.OSS} is not satisfied, and hence this theorem  gives  no information on whether \(E\cdot E\) is a \(U\)-set or an \(M_0\)-set. This together with Theorem \ref{rem:bounded-all-dimensions} indicates that the energy condition in Theorem \ref{thm.OSS} is only a sufficient condition for the arithmetic product set to be an $M_0$-set.

\item Since the $U$-set $E$ constructed in Theorem \ref{thm:sharpness}(a) has Hausdorff dimension $1/2$ and $E\cdot E$ is an $H$ set,   it follows that the main contribution    of  Theorem \ref{thm:sharpness}(b) is to give a stronger bound 	$\underline\dim_BE\le1/2$ under the hypothesis that  the self-product $E\cdot E$ is an $H$-set. Moreover, the $H$-set condition in \ref{thm:sharpness}(b) can be replaced by an $H^{(n)}$-set condition to yield an estimate for the lower box dimension; see Proposition \ref{prop:hn-bound}. In this setting, \eqref{eq:hn-bound} gives the sharp endpoint   $1/2$ provided that the relevant frequency windows satisfy  	$\liminf_{j\to\infty}\log U_j/\log L_j=1$.
	\end{itemize}
  	
	The rest of this paper is organized as follows. Section \ref{sec:preliminaries} collects the product criterion, a criterion for the set of uniqueness and dimension formulas for homogeneous  Moran sets. In Section  \ref{sec:self-similar-u-sets}, we prove theorem \ref{thm:intro-self-similar-product}; In Section \ref{sec:bounded-integer-moran}, we prove theorems \ref{thm:intro-bounded-integer-moran} and \ref{rem:bounded-all-dimensions}. In the last two  sections, we will prove Theorems   \ref{thm:unbounded-m0-product}  and \ref{thm:sharpness}, respectively.

\subsection{Notation} Throughout this paper, we use the following standard  notation: for any real number $x\in\mathbb R$, $\lfloor x\rfloor$  denotes the floor function, i.e., the greatest integer less than or equal to $x$; $\lceil x\rceil$ denotes the ceiling function, i.e., the least integer greater than or equal to $x$; \(\{x\}\) denotes the fractional part of $x$, i.e., $\{x\}=x-\lfloor x\rfloor;$ \(\|x\|\) denotes the distance from $x$ to the nearest integer, i.e., $\|x\|=\min_{n\in\mathbb Z} |x-n|$. Equivalently, $\|x\|=\min\{\{x\},1-\{x\}\},$ so that $0\le \|x\|\le \frac12.$
	
	\section{Preliminaries}	\label{sec:preliminaries}
	In this section, we present some basic preliminaries required for the proofs of our main results.

	\subsection{A criterion for $M_0$-set}
	The following lemma extends  \cite[Lemma 2]{Kaufman1993} to arbitrary Borel probability measures, rather than the middle-third Cantor measure.
	\begin{lemma}\label{lem:product-kaufman}
		Let \(E\subset\mathbb R\) be a compact set.  Suppose that there exist a probability measure \(\mu\) supported on \(E\) and a countable set \(\Lambda\subset\mathbb R\)   such that
		\begin{equation}\label{eq:product-cross-decay}
	\lim_{|u|\to\infty}\widehat\mu(u)\,\widehat\mu(\tau u)=0
	\qquad \text{for every } \tau\notin\Lambda.
		\end{equation}
		Thus, for every uncountable compact set \(F\subset\mathbb R\), the set \(F\cdot E\) is an \(M_0\)-set.
	\end{lemma}
	
	\begin{proof}
		If \(F\) is an uncountable compact set, then $F$ supports a continuous Borel probability measure $\sigma$ (see, e.g., \cite[p. 53, Theorem 8.1]{Par1967}). Let \(\rho(\cdot)=(\sigma\times\mu)(h^{-1}(\cdot))\), where $h(s,t)=st$ with $s,t\in\mathbb{R}$, be the Borel probability measure supported on $F\cdot E$, i.e.,
\[\int \varphi(x)\,d\rho(x)=\iint \varphi(st)\,d\sigma(s)d\mu(t),	\qquad 	\varphi\in C_c(\mathbb R).\]
Thus, the Fourier transform of $\rho$ satisfies that
\[\widehat\rho(u)=\int e^{-2\pi i ux}d\rho(x)=\iint e^{-2\pi iust}\,d\sigma(s)d\mu(t)=\int \widehat\mu(su)\,d\sigma(s),\] which gives
		\begin{equation}\label{eq:rho-cauchy-bound}
			|\widehat\rho(u)|^2 \le \iint |\widehat\mu(su)\,\widehat\mu(tu)|\,d\sigma(s)d\sigma(t).
		\end{equation}	Let\[\Gamma=\{(s,t)\in\mathbb{R}^2:t=0\}\bigcup\bigcup_{\alpha\in\Lambda}\{(s,t)\in\mathbb{R}^2:s=\alpha t\}=:E^0 \bigcup\bigcup_{\alpha\in\Lambda} E_\alpha.\]	Since \(\sigma\) has no atoms, it follows from Fubini's theorem that \[(\sigma\times\sigma)(E_\alpha)=\int \sigma(\{\alpha t\})\,d\sigma(t)=0\] for all $\alpha\in\Lambda$, and the same holds for $E^0$ as well. Since \(\Lambda\) is countable,  we have  \((\sigma\times\sigma)(\Gamma)=0\). Therefore, if \((s,t)\notin \Gamma\), then \(t\neq0\) and \(s/t\notin\Lambda\), and  thus the condition \eqref{eq:product-cross-decay} yields
		\[\widehat\mu(su)\,\widehat\mu(tu)=\widehat\mu\left(\frac{s}{t}tu\right)\widehat\mu(tu)\to0 \qquad (|u|\to\infty).\]
		Note that \(|\widehat\mu(su)\,\widehat\mu(tu)|\le1\) for all $u\in\mathbb{R}$. Applying the dominated convergence theorem to \eqref{eq:rho-cauchy-bound}, we obtain that \(\widehat\rho(u)\to0\) if ${|u|\to\infty}$. Thus, \(\rho\) is a Rajchman measure supported on \(F\cdot E\), and hence \(F\cdot E\) is an $M_0$-set. This finishes the proof.
	\end{proof}
	
	\begin{remark}\label{rem:exceptional}
		Let us briefly explain the hypothesis and applications of Lemma \ref{lem:product-kaufman}.
		\begin{enumerate}
			\item[(1)] By the properties of the Fourier transform of measures, the condition  \eqref{eq:product-cross-decay} actually says that  if \(\tau\notin\Lambda\), then the convolution measure $\mu\ast (\mu \circ S_\tau^{-1})$ is a Rajchman measure supported on the set $E+\tau E$, where $S_\tau(x)=\tau x$ is a scaling map and $E$ is the support of the measure $\mu$.
			\item[(2)] In the proofs of Theorem \ref{thm:intro-self-similar-product} and Theorem \ref{thm:intro-bounded-integer-moran}, one of the crucial steps is to choose a suitable measure $\mu$ so that \eqref{eq:product-cross-decay} holds outside a countable exceptional set $\Lambda$, where $\Lambda$ is taken to be $\mathbb Q(\lambda)$ and \(\mathbb Q\), respectively.
		\end{enumerate}
	\end{remark}
	
	\subsection{Criteria  for $U$-set}
	We recall a classical criterion of Rajchman.
	\begin{theorem}[Rajchman {\cite[p. 50, Theorem II]{Salem1963}}; {\cite[Chapter IX, Theorem 6.3]{Zyg1959}}]\label{thm:salem-zygmund-criterion}
		If the compact set  \(E\subset[0,1]\) is an $H$-set, i.e., there exist a strictly increasing sequence of positive integers \(n_j\to\infty\) and a nonempty open interval \(I\subset(0,1)\) such that
		\[\{n_jx\}\notin I\qquad (x\in E,\ j\ge1),\]
		then \(E\) is a \(U\)-set.
	\end{theorem}
	The following result is well-known (see, e.g., \cite{LaiXie2026}), we provide the proof for completeness. Here the integers $m_k\ge2$ need not be bounded.
	\begin{lemma}\label{lem:left-half-moran-u}
		Suppose that \(\frac{\max \{d\in D_k\}}{m_k}\le L<1/2\) for every \(k\ge1\). Then \(K\) is a \(U\)-set.
	\end{lemma}
	\begin{proof}
		Put \(Q_k=m_1\cdots m_k\) for each $k$.  If \(x=\sum_{k=1}^{\infty}\frac{d_k}{Q_k}\in K,\) then
		\[Q_jx=\sum_{k=1}^{j}d_k{m_{k+1}\cdots m_j}+\sum_{k=j+1}^{\infty}\frac{d_k}{m_{j+1}\cdots m_k}.\]
		The first sum is an integer.  Hence
\[\{Q_jx\}=\sum_{k=j+1}^{\infty}\frac{d_k}{m_{j+1}\cdots m_k}\le L\left(1+\frac1{m_{j+1}}+\frac1{m_{j+1}m_{j+2}}+\cdots\right)		\le 2L.\]
Since \(2L<1\), we have \(\{Q_jx\}\notin(2L,1)\) for every \(x\in K\) and \(j\ge1\). 	Theorem \ref{thm:salem-zygmund-criterion}, applied with
		\(n_j=Q_j\) and \(I=(2L,1)\), implies that \(K\)	is an $H$-set, and hence a \(U\)-set.
	\end{proof}
	
\subsection{Dimension formulas for homogeneous Moran sets}
Let \(E \subset \mathbb{R}\) be a nonempty Borel set. The diameter of $E$ is defined as $|E|=\sup\{|x-y|: x, y \in E\}$, that is, the greatest distance apart of any pair of points in $E.$ For \(s \ge 0\) and \(0 < \delta \le 1\), the \(s\)-Hausdorff measure of \(E\) is defined by
\[\mathcal{H}^{s}(E) = \lim_{\delta \to 0} \mathcal{H}_{\delta}^{s}(E),\]where \[\mathcal{H}_{\delta}^{s}(E) = \inf \left\{ \sum_{i=1}^{\infty} |U_i|^s : E\subset \bigcup_{i} U_i, \ |U_i| < \delta \right\}.\]
The Hausdorff dimension of \(E\) is
$$\dim_H E = \inf\{s \ge 0 : \mathcal{H}^{s}(E) = 0\} = \sup\{s \ge 0 : \mathcal{H}^{s}(E) = \infty\}.$$

	The formulas in \textup{(i)} and \textup{(ii)} are consequences of \cite[Theorem 2.1 and Corollary 2.1]{FengWenWu1997} and \cite[Lemma 2.2]{FengWenWu1997}, respectively.
	\begin{lemma}[Hausdorff dimension formulas]\label{lem:moran-dimension}
Let $K$ be the homogeneous 	Moran set as in \eqref{1.12}, and write $N_n=\#D_n$, $A_n=\prod_{k=1}^nN_k$, and $Q_n=\prod_{k=1}^nm_k$.		\begin{enumerate}[label=\textup{(\roman*)}]
\item If $\sup_n m_n<\infty$, then
\[\dim_H K=\liminf_{n\to\infty}\frac{\log A_n}{\log Q_n}.\]			\item If $D_n=\{0,1,\ldots,N_n-1\}$ for every $n$, then	
\[\dim_H K=\liminf_{n\to\infty}\frac{\log A_n}{-\log \frac{N_{n+1}}{Q_{n+1}}}.\]
\end{enumerate}	\end{lemma}

	\section{On self-similar sets of uniqueness}
	\label{sec:self-similar-u-sets}
	The goal of this section is to prove Theorem \ref{thm:intro-self-similar-product}.  Our argument depends on Lemma \ref{lem:product-kaufman} and  a theorem of   Varj\'u and Yu \cite[Theorem 1.5]{VarjuYu2022} on the Fourier decay of self-similar measures.
	
	\noindent\textit{Principal Assumptions}.  In what follows,  for generality, we shall usually suppose that  $k\geq2$ is a positive integer and  \[\{f_i(x)=r^{\ell_i}x+a_i:1\le i\le k\}\] is an IFS satisfying   \(0<r<1\), \(\ell_i\in\mathbb Z_{>0}\) with  \(\gcd(\ell_1,\dots,\ell_k)=1\), and $a_1,\ldots,a_k\in\mathbb{R}.$ Let  $E$ be the non-singleton self-similar set of this IFS. Let \((p_1,\ldots,p_k)\) be a probability vector with strictly positive entries and let \(\mu\) be the corresponding self-similar measure, i.e., \(\mu=\sum_{i=1}^k p_i\mu\circ f_i^{-1}\).

For each finite word \(\mathbf i=(i_1,\ldots,i_s)\in\{1,2,\ldots,k\}^s\), we define the cylinder mappings
	\[f_{\mathbf i}(x)=\bigl(f_{i_1}\circ\cdots\circ f_{i_s}\bigr)(x)
	=r^{L(\mathbf i)}x+a_{\mathbf i},\]
	where \(L(\mathbf i)=\ell_{i_1}+\cdots+\ell_{i_s}\) is the contraction exponent and
	\(a_{\mathbf i}=a_{i_1}+r^{\ell_{i_1}}a_{i_2}+\cdots+r^{\ell_{i_1}+\cdots+\ell_{i_{s-1}}}a_{i_s}\)
	is the cylinder translation.
	
	We first provide the following elementary separation lemma for cylinder mappings.
	\begin{lemma}\label{lem:separation}
		If the self-similar set \(E\) is not a singleton, then there exist two finite words \(\mathbf u,\mathbf v\) and an integer \(N\ge1\) such that \[L(\mathbf u)=L(\mathbf v)=N \quad \text{and} \quad a_{\mathbf u}\ne a_{\mathbf v}.\]
	\end{lemma}
	
	\begin{proof}
		Suppose, to the contrary, that  \(a_{\mathbf u}=a_{\mathbf v}\) whenever \(L(\mathbf u)=L(\mathbf v)\).  Pick \(i\ne j\in\{1,2,\ldots,k\}\).  Let \(\mathbf u\) be the word consisting of \(\ell_j\) copies of \(i\), and let 	\(\mathbf v\) be the word consisting of \(\ell_i\) copies of \(j\), i.e.,
		$$\mathbf u=(\underset{\ell_j ~ {\rm{times}}}{\underbrace{i,\ldots,i}}), \qquad
		\mathbf v=(\underset{\ell_i~ {\rm{times}}}{\underbrace{j,\ldots,j}}).$$
		Then \(L(\mathbf u)=L(\mathbf v)=\ell_i\ell_j\), and hence \(a_{\mathbf u}=a_{\mathbf v}\). However, a simple calculation gives
		\begin{align*}
a_{\mathbf u}&=a_i\frac{1-r^{\ell_i\ell_j}}{1-r^{\ell_i}}\qquad \text{and} \qquad		
a_{\mathbf v}=a_j\frac{1-r^{\ell_i\ell_j}}{1-r^{\ell_j}}.
		\end{align*}
		It follows from  \(0<r<1\) that \[\frac{a_i}{1-r^{\ell_i}}=\frac{a_j}{1-r^{\ell_j}},\] which is exactly the fixed point of the maps $f_i$ and $f_j$. Since the distinct indices \(i,j\) are  arbitrary, the set $E$ is therefore a singleton, a contradiction. This finishes the proof of Lemma \ref{lem:separation}.
	\end{proof}

The proof next is also based on a theorem of Varj\'u and Yu \cite{VarjuYu2022}, which gives a   characterization on Fourier decay of self-similar measures.  Let $\exp(x):=e^x$.
	
	\begin{theorem}[Varj\'u and Yu {\cite[Theorem 1.5]{VarjuYu2022}}]
		\label{thm:vy-scale-uniform}
		Suppose the Principal Assumptions hold,  \(\ell_1=\ell_2\) and  \(\Delta:=a_1-a_2\not=0\). Then there exists  a constant	\(C>0\) depending only on \(\ell_1,\dots,\ell_k\) and \(p_1,\ldots,p_k\) such that \[|\widehat\mu(u)|\le \exp\left(-C^{-1}\sum_{j>C}\|\Delta u r^j\|^2\right).\]
	\end{theorem}
	\begin{remark}
		If, in addition, \(\lambda=r^{-1}\) is a Pisot number, \(a_1,\ldots,a_k\in\mathbb Q(\lambda)\), and the associated self-similar set \(E\) has zero Lebesgue measure, by Theorem \ref{thm:intro-varju-yu-form}, the set $E$ is a $U$-set, so every self-similar measure $\mu$ is not a Rajchman measure, i.e., $\widehat{\mu}(u)\not\to0$ as $|u|\to\infty.$  Nevertheless, our key Lemma  \ref{lem:self-similar-cross} below shows that for each self-similar measure $\mu$ in Theorem \ref{thm:vy-scale-uniform},  if  \(\tau\notin\mathbb Q(\lambda)\), then $\widehat\mu(u)\widehat\mu(\tau u)\to0$ as $|u|\to\infty.$

	\end{remark}
	According to Lemma \ref{lem:separation} and Theorem \ref{thm:vy-scale-uniform}, we get a uniform estimate of the following series of squared distances.
	\begin{proposition}\label{prop:pisot-case-vy-estimate}
		Suppose the Principal Assumptions hold. There exist a positive number  \(b\), depending on the fixed IFS, and an integer \(J\ge0\),  depending only on \(\ell_1,\dots,\ell_k\) and \(p_1,\ldots,p_k\), such that for each \(\delta>0\) there is a constant \(C_\delta<\infty\), depending only on $\delta$, such that
		\begin{equation}
			\label{eq:vy-uniform-estimate}
			|\widehat\mu(u)|\ge\delta \quad \implies	\quad \sum_{j\ge J}\|bu r^{j}\|^2\le C_\delta.
		\end{equation}
	\end{proposition}
	
	\begin{proof}
		Applying Lemma \ref{lem:separation} to the IFS \(\{f_i(x)=r^{\ell_i}x+a_i:1\le i\le k\}\), there are two words \(\mathbf i\ne\mathbf j\) and a positive integer \(N\ge1\) such that	\(L(\mathbf i)=L(\mathbf j)=N\) and \(a_{\mathbf i}>a_{\mathbf j}\).
		
		Put \(b_0=a_{\mathbf i}-a_{\mathbf j}>0\), and set \(\mathbf w_1=\mathbf i\mathbf j\),
		\(\mathbf w_2=\mathbf j\mathbf i\). Note that the words $\mathbf w_1, \mathbf w_2$ have the same length, we denote it by \(|\mathbf w_1|=|\mathbf w_2|=m\) for some $m$. Therefore,
		\[L(\mathbf w_1)=L(\mathbf w_2)=2N, \qquad a_{\mathbf w_1}=a_{\mathbf i}+r^Na_{\mathbf j}, \quad a_{\mathbf w_2}=a_{\mathbf j}+r^Na_{\mathbf i}.\]
		Let \[\widetilde b:=a_{\mathbf w_1}-a_{\mathbf w_2}=(1-r^N)b_0>0.\]

		Iterating the self-similar identity \(\mu=\sum_{i=1}^k p_i\mu\circ f_i^{-1}\)   gives
		\[\mu=\sum_{|\mathbf w|=m}p_{\mathbf w}\mu\circ f_{\mathbf w}^{-1}
		\quad \text{where } \quad	p_{\mathbf w}=p_{w_1}\cdots p_{w_m}>0.\]
		Thus, the self-similar measure $\mu$ can be represented as
		\[\mu=\frac12\sum_{i=1}^k p_i\mu\circ f_i^{-1}+\frac12\sum_{|\mathbf w|=m}p_{\mathbf w}\mu\circ f_{\mathbf w}^{-1},\] which corresponds  to the following IFS
		$$\{f_i(x)=r^{\ell_i}x+a_i:1\le i\le k\}\cup \{f_{\mathbf w}(x)=r^{L({\mathbf w})}x+a_{\mathbf w}: |{\mathbf w}|=m, {\mathbf w}\in \{1,\ldots,k\}^m\},$$
		which contains, in particular, the two mappings
		\[
		f_{\mathbf w_1}(x)=r^{2N}x+a_{\mathbf w_1},
		\qquad
		f_{\mathbf w_2}(x)=r^{2N}x+a_{\mathbf w_2},
		\qquad
		a_{\mathbf w_1}-a_{\mathbf w_2}=\widetilde b.
		\]
		
		Because their contraction exponents satisfy
		\[
		\gcd\Bigl(\{\ell_1,\ldots,\ell_k\}\cup\{L(\mathbf w):|\mathbf w|=m\}\Bigr)
		=\gcd(\ell_1,\ldots,\ell_k)=1,
		\]
		by relabelling the exponent if necessary, it follows from Theorem \ref{thm:vy-scale-uniform} that there is a constant \(C>0\), depending only on $\ell_1,\ldots,\ell_k, L(\mathbf w)$ and $p_1,\ldots,p_k, p_{\mathbf w}$ with $|\mathbf w|=m$, such that
		\[|\widehat\mu(u)|\le \exp\left(-C^{-1} \sum_{j>C}\|\widetilde b\,u r^j\|^2\right).\]
		The definitions of $L(\mathbf w)$ and $p_{\mathbf w}$ imply that $C$ depends only on \(\ell_1,\dots,\ell_k\) and \(p_1,\ldots,p_k\). Hence  \(|\widehat\mu(u)|\ge\delta\) yields \[\sum_{j>C}\|\widetilde b\,u r^{j}\|^2\le -C\log \delta.\]
		Taking \(J=\lfloor C\rfloor+1\), \(b=\widetilde b\), and \(C_\delta= -C\log \delta\) proves \eqref{eq:vy-uniform-estimate}.  The proof of Proposition \ref{prop:pisot-case-vy-estimate} is complete.
	\end{proof}

	 lemma \ref{lem:self-similar-cross} shows that the condition \eqref{eq:product-cross-decay} of Lemma \ref{lem:product-kaufman} holds for self-similar measures outside the exceptional set $\mathbb Q(\lambda)$.
	\begin{lemma}\label{lem:self-similar-cross}
		Suppose the Principal Assumptions hold. Thus, \(\lambda:=r^{-1}\) is a Pisot number, and for every \(\tau\notin\mathbb Q(\lambda)\),
		\[\widehat\mu(u)\widehat\mu(\tau u)\to0 \qquad (|u|\to\infty).\]
	\end{lemma}
	
	The proof of Lemma \ref{lem:self-similar-cross} relies on the following classical rigidity theorem of Pisot.
	\begin{lemma}[Pisot {\cite[Chapter 1, Theorem A]{Salem1963}}]\label{lem:pisot-rigidity}
		If \(\lambda>1\) is such that there is a \(v\in\mathbb R\setminus\{0\}\) such that
		\[\sum_{n=0}^{\infty}\|\lambda^n v\|^2<\infty,\]
		then \(\lambda\) is a Pisot number and \(v\in\mathbb Q(\lambda).\)
	\end{lemma}
	
	\begin{proof}[Proof of Lemma \ref{lem:self-similar-cross}]
		Since $|\widehat{\mu}(u)|=|\widehat{\mu}(-u)|$ for all $u\in\mathbb{R}$, we may suppose on the contrary that there is a positive number \(\tau\notin\mathbb Q(\lambda)\), a constant \(\delta>0\)  and a sequence
		\(u_n\to+\infty\) such that
		\[ |\widehat\mu(u_n)|\ge\delta, \qquad |\widehat\mu(\tau u_n)|\ge\delta.\]
		Let \(b,J\) be as in Proposition \ref{prop:pisot-case-vy-estimate}, and $C_\delta$ be the constant such that
		\begin{align}\label{eq:self-similar-square-bound}
			\sum_{j\ge J}\|bu_n\lambda^{-j}\|^2 &\le C_\delta \qquad \text {and} \qquad
			\sum_{j\ge J}\|b\tau u_n\lambda^{-j}\|^2\le C_\delta.
		\end{align}
		Let \(N_n,M_n\) be the natural integers increasing to $+\infty$ as $n\to+\infty$ such that
		\[v_n:=bu_n\lambda^{-N_n}\in[1,\lambda), \qquad w_n:=b\tau u_n\lambda^{-M_n}\in[1,\lambda).\]
		
		Keeping a positive integer \(M\ge0\) fixed and letting \(n\) sufficiently large such that
		\(N_n-M\ge J\) and  \(M_n-M\ge J\), so \eqref{eq:self-similar-square-bound} gives
		\[\sum_{q=0}^{M}\|v_n\lambda^q\|^2=\sum_{q=0}^{M}\|bu_n\lambda^{q-N_n}\|^2\le C_\delta\quad \text{and} \quad \sum_{q=0}^{M}\|w_n\lambda^q\|^2=\sum_{q=0}^{M}\|b\tau u_n\lambda^{q-M_n}\|^2\le C_\delta.\]
		
		By extracting, if necessary, subsequences from $\{v_n\}$ and $\{w_n\}$ respectively,
		we may assume that $v_n\to v\in[1,\lambda]$ and $w_n\to w\in[1,\lambda]$.
		Letting $n\to\infty$ and then $M\to\infty$, we have
		\[\sum_{q=0}^{\infty}\|v\lambda^q\|^2<\infty\quad \text{and} \quad  	\sum_{q=0}^{\infty}\|w\lambda^q\|^2<\infty. \]
		By Lemma \ref{lem:pisot-rigidity},  $\lambda$ is a Pisot number, and  \(v,w\in\mathbb Q(\lambda)\).
		
		Finally, note that \(v_n,w_n\in[1,\lambda)\) implies that $\lambda^{-1}<\frac{w_n}{v_n}=\tau\lambda^{N_n-M_n}<\lambda$, it follows that  $$-(1+\log_\lambda {\tau})<N_n-M_n<1-\log_\lambda {\tau}.$$ Passing to a further subsequence, if necessary, we may assume that \(N_n-M_n=q_0\) for some integer  $q_0$ and for all $n$.  Letting \(n\to\infty\) in $\frac{w_n}{v_n}=\tau\lambda^{N_n-M_n}$, we obtain
		\[\tau=\lambda^{-q_0}\frac{w}{v}\in\mathbb Q(\lambda).\]
		It is a contradiction to \(\tau\notin\mathbb Q(\lambda)\). This completes the proof of Lemma \ref{lem:self-similar-cross}.
	\end{proof}
	
	Now we can prove Theorem \ref{thm:intro-self-similar-product}.
	\begin{proof}[Proof of Theorem \ref{thm:intro-self-similar-product}]
		It suffices to prove the sufficiency.	Let \(E\subset\mathbb R\) be a non-singleton self-similar set of uniqueness, which is generated by the IFS $\{g_i(x)=r_i x+b_i:1\le i\le k, r_i\in (0,1), b_i\in\R\}$. Thus, $E=\bigcup_{i=1}^{k}g_i(E)$.  By Theorem \ref{thm:intro-varju-yu-form},  the IFS  $\{g_i\}_{i=1}^k$ is conjugated to the following Pisot form IFS:
		\[f_i(x)=r^{\ell_i}x+a_i,\qquad i=1,\dots,k,\] where \(0<r<1\), \(\ell_i\in\mathbb Z_{>0}\), 	\(\gcd(\ell_1,\dots,\ell_k)=1\), \(\lambda=r^{-1}\) is a Pisot number, and \(a_i\in\mathbb Q(\lambda)\) for all \(i\). That is to say, there is a similarity \(S(x)=\alpha x+\beta\), \(\alpha\ne0\), such that for each $i=1,\ldots,k$,
		\begin{equation}\label{eq.equv.}
			S\circ g_i\circ S^{-1}=f_i \quad \text{or} \quad g_i\circ S^{-1}=S^{-1} \circ f_i\quad \text{or} \quad  g_i^{-1}\circ S^{-1}=S^{-1}\circ f_i^{-1}.
		\end{equation}
		Let $F$  be the self-similar set generated by the IFS $\{f_i\}_{i=1}^k,$ that is, $F=\bigcup_{i=1}^{k}f_i(F)$. By \eqref{eq.equv.}, we have
		$$\cup_{i=1}^k g_i(S^{-1}(F))=\cup_{i=1}^k S^{-1}(f_i(F))=S^{-1}\cup_{i=1}^k (f_i(F))=S^{-1}(F).$$ The uniqueness of the self-similar set of $\{g_i\}_{i=1}^k$ implies that $E=S^{-1}(F)$. Since two homothetic sets (i.e., one is a scaled and translated copy of the other) share the same $U$-set property (see, e.g., \cite[p. 61, Theorem VII]{KahaneSalem1963}), $F$ is a self-similar set of uniqueness.

		Given a positive probability vector $(p_1,\ldots,p_k),$ let $\mu_0$ and $\mu$ be the associated self-similar measure supported on the set $E$ and $F$, respectively. Then $\mu_0=\sum_{i=1}^k p_i\mu_0\circ g_i^{-1}$ and $\mu=\sum_{i=1}^kp_i\mu \circ f_i^{-1}$. Therefore, by \eqref{eq.equv.}, for any Borel set $A\subset\R$,
		\begin{equation*} \mu_0(S^{-1}A)=\sum_{i=1}^kp_i\mu_0(g_i^{-1}S^{-1}A)=\sum_{i=1}^kp_i\mu_0(S^{-1}f_i^{-1}A),
		\end{equation*}
		which yields that $\mu_0\circ S^{-1}=\mu$ by the uniqueness of the self-similar measure of the  IFS $\{f_i\}_{i=1}^k$.  Applying Lemma \ref{lem:self-similar-cross} to the  IFS $\{f_i(x)\}_{i=1}^k$, we obtain that, for every \(\tau\notin\mathbb Q(\lambda)\),
		\[
		\widehat{\mu}(u)\widehat{\mu}(\tau u)\to0
		\qquad (|u|\to\infty).
		\]
		Note that the Fourier transform of $\mu$ satisfies that
		$$\widehat{\mu}(u)=\int e^{-2\pi i ut}d\mu(t)=e^{-2\pi iu\beta} \widehat{\mu_0}(\alpha u).$$
		Therefore,  for every \(\tau\notin\mathbb Q(\lambda)\),	
		\[\widehat{\mu_0}(u)\widehat{\mu_0}(\tau u)\to0		\qquad (|u|\to\infty),\]
		hence the sufficiency follows from Lemma \ref{lem:product-kaufman}, as required. The proof of Theorem \ref{thm:intro-self-similar-product} is complete.
	\end{proof}

	\begin{remark}
		For the non-singleton self-similar set $E$ in Principal Assumptions, from the proof of Theorem \ref{thm:intro-self-similar-product}, we know that,   		without the uniqueness assumption on $E$,  \[F\text{ is uncountable}\quad\Longrightarrow\quad F\cdot E \text{ is an  $M_0$-set}.\]
		However, Example \ref{ex:pisot-interval} below shows that the uniqueness assumption is essential for the implication \[F\cdot E \text{ is an  $M_0$-set} \quad\Longrightarrow\quad F\text{ is uncountable}.\]
	\end{remark}
	
	At the end of this section, we apply Theorem \ref{thm:intro-self-similar-product}  to the homogeneous IFS \(\{\phi_i(x)=r x+ a_i:1\le i\le m\}\), where \(0<r<1\) and  \(A=\{a_1,\dots,a_m\}\subset\mathbb R\) with \(m=\#A\ge2\). Then the corresponding homogeneous self-similar set \(C_{r,A}=\bigcup_{i=1}^{m}\phi_i(C_{r,A})\) can be represented as  \[C_{r,A}=\left\{	\sum_{n=0}^{\infty}a_{i_n}r^{n}:i_n\in\{1,\dots,m\}	\right\}.\]
	In particular, the \emph{middle-$r$ Cantor set} \(C_{r}\) mentioned in Subsection \ref{sec.1} is obtained if we take $A=\{0,1-r\}$. Recall that \(C_r\) is a set of uniqueness if and only if $r^{-1}$ is a Pisot number by Salem and Zygmund's theorem. The following result follows immediately from Theorem \ref{thm:intro-self-similar-product}.
	\begin{theorem}\label{thm:intro-pisot-cantor-product}
		If \(C_{r,A}\) is a set of uniqueness, then for every compact set	\(F\subset\mathbb R\),
		\[F\cdot C_{r,A} \text{ is an  $M_0$-set} \quad\text{if and only if}\quad F\text{ is uncountable}.\]
	\end{theorem}

	\begin{remark}\label{rem:checking-U-property}
		For the homogeneous self-similar set, if we assume additionally that $\lambda=r^{-1}$ is a Pisot number and $A\subset\mathbb Q(\lambda)$, since \(\#A\ge2\),  by Theorem \ref{thm:intro-varju-yu-form} of Varj\'u and Yu \cite{VarjuYu2022},  we get that \(C_{r,A}\) is a \(U\)-set if and only if $C_{r,A}$ has  zero Lebesgue measure, i.e., \(\mathcal L^1(C_{r,A})=0.\) So the uniqueness assumption in Theorem \ref{thm:intro-pisot-cantor-product} can be	replaced by \(\mathcal L^1(C_{r,A})=0.\)
	\end{remark}

	\begin{example}[Necessity of the uniqueness assumption]
		\label{ex:pisot-interval}
		Let \(\lambda=r^{-1}=\frac{1+\sqrt5}{2}, A=\{0,r\}.\) Then \(\lambda\) is a non-integer Pisot number and \(A\subset\mathbb Q(\lambda)\).  Denote the corresponding IFS by
		\[\phi_0(x)=r x, \qquad \phi_1(x)=r (x+1).\] Obviously, $\phi_0([0,\lambda])=[0,1]$.
		Since \(1+\lambda^{-1}=\lambda\), it follows that 	\(\phi_1([0,\lambda])=[r,\lambda].\)
		Therefore, $[0,\lambda]=\cup_{i=0}^1 \phi_i([0,\lambda])$.  In this case, the self-similar set \(C_{r,A}:=[0,\lambda]\) is an \(M_0\)-set.  Taking \(F=\{1\}\), we obtain \(F\cdot C_{r,A}=C_{r,A}\) is an \(M_0\)-set, but \(F\) is finite.  Therefore the implication
		\[F\cdot C_{r,A}\text{ is an $M_0$-set}\quad \Longrightarrow\quad F\text{ is uncountable}\]
		fails without the assumption that \(C_{r,A}\) is a \(U\)-set.
	\end{example}

	\section{On bounded integer-based homogeneous Moran sets}
	\label{sec:bounded-integer-moran}
	The goal of this section is to prove Theorem \ref{thm:intro-bounded-integer-moran} and Theorem  \ref{rem:bounded-all-dimensions}. In this section, we always assume that  the homogeneous Moran set \(K\) defined as in \eqref{1.12} satisfies the assumptions of Theorem \ref{thm:intro-bounded-integer-moran}, and let  $\mu:=\overset{\infty}{\underset{k=1}{\scalebox{2}{\(\ast\)}}} \left(\frac{1}{\#D_k}\sum_{d_k\in D_k} \delta_{d_ky_k}\right)$ (as in \eqref{eq.4.1}) be the natural probability measure supported on $K$, where  \(y_k=Q_k^{-1}\) and \(Q_k=m_1\cdots m_k\). It is not difficult to compute that the Fourier transform of the measure $\mu$ is 	\[\widehat\mu(u)=\prod_{k=1}^{\infty}P_{D_k}(y_k u),\]	where the symbol $P_D$ denotes the mask polynomial  $P_D(u)=\frac1{\#D}\sum_{d\in D}e^{-2\pi i d u}.$
	
The lemma below gives the resulting uniform bound of the functions $P_D$.
	\begin{lemma}\label{lem:uniform-mask-estimate}
		There exists a constant \(c>0\), depending only on \(B\), such that for every integer \(m\) with \(2\le m\le B\) and every digit set \(D\subset\{0,1,\ldots,m-1\}\) with \(\#D\ge2\),
		\[|P_D(u)|\le \exp\left(-c\|g(D)u\|^{2}\right),\qquad (u\in\mathbb R),\]
		where \(g(D)\) denotes the positive greatest common divisor of the differences \(d-d'\), \(d\not=d'\in D\).
	\end{lemma}
	
	\begin{proof}
		Let $g(D):=\gcd\{d-d'>0: d,d'\in D, d\not=d'\}$, \(d_0=\min\{d\in D\}\), and put
		\[D'=\left\{\frac{d-d_0}{g(D)}:d\in D\right\}\subset\mathbb{Z}_{\geq0}.\]
		Then  \(\gcd(D')=1\) and \(D=d_0+g(D)D'\). Therefore,
		\[P_D(u)=e^{-2\pi i d_0u}P_{D'}(g(D)u),\]
		hence it suffices to estimate the upper bound of the function  $|P_{D'}(u)|$.
		
		First, using the Taylor expansion of the exponential function, if $u\to0$, then
		\[\begin{aligned}
|P_{D'}(u)|^2&=\frac1{(\#D')^2}\sum_{a,b\in D'} e^{-2\pi i(a-b)u}
=1-\frac{2\pi^2}{(\#D')^2}\sum_{a,b\in D'}(a-b)^2u^2+O(u^4).
		\end{aligned}\]
		Notice that $C:=\frac{2\pi^2}{(\#D')^2}\sum_{a,b\in D'}(a-b)^2>0$ since $\#D'\geq2,$ and $P_{D'}$ is an integer-periodic function.  Thus, there are constants $c_1\in (0, C)$ and $\eta\in (0,1/2)$ such that
		$$|P_{D'}(u)|^2 \leq \exp\left({-c_1\|u\|^2}\right) \qquad \text{for all } \|u\|<\eta.$$
		Let $c_{D'}=\frac{c_1}{2}$. We obtain  $$|P_{D'}(u)| \leq   \exp\left({-c_{D'}\|u\|^2}\right) \qquad \text{for all } \|u\|<\eta.$$
		
		Second, it is not hard to see that  \(\gcd(D')=1\) implies that \(|P_{D'}(u)|=1\) if and only if \(u\in\mathbb Z\). This means that  $|P_{D'}(u)|<1$ for all $\eta\leq\|u\|\leq\frac{1}{2}$.
		Whence, decreasing $c_{D'}$ if necessary, we get that
		\[|P_{D'}(u)|\leq \exp\left(-c_{D'}\|u\|^{2}\right), \qquad (u\in\mathbb R).\]
		
		Finally, since \(D\subset\{0,1,\ldots,m-1\}\) and \(m\le B\), there are only finitely many such digit set \(D'\).  Taking the minimum of the corresponding positive constants $c_{D'}$ gives a constant \(c>0\) depending only on \(B\), as desired. The proof of Lemma \ref{lem:uniform-mask-estimate} is complete.
	\end{proof}
	
	\begin{lemma}\label{lem:multidigit-integer-rigidity}
		Fix a positive number $\delta>0$, suppose that \(m_1\leq |u_j|\to\infty\) and
		\(|\widehat\mu(u_j)|\ge\delta>0\). Thus, for each integer \(j\), there is an index \(p_j\) such that \(1\le |y_{p_j}u_j|\leq B\), and every accumulation point of \(\{y_{p_j}u_j\}\) is a rational number.
	\end{lemma}
	
	\begin{proof}
		As the sequence $\{m_1\ldots m_k\}_{k=1}^\infty$ is strictly increasing as $k\to\infty$ and $|u_j|\geq m_1$ for all $j$, there is a unique integer $p_j$  such that $m_1\ldots m_{p_j}\leq |u_j|< m_1\ldots m_{p_j+1}$. Thus, \(1\le |y_{p_j}u_j|\le m_{p_j+1}\leq B\).

		On the one hand,  it follows from Lemma \ref{lem:uniform-mask-estimate} that
		\[|\widehat\mu(u_j)|=\prod_{k=1}^{\infty}|P_{D_k}(y_k u_j)|
		\le \exp\left(-c\sum_{k=1}^{\infty}	\|g(D_k)y_k u_j\|^{2}\right).\]
		By \(|\widehat\mu(u_j)|\ge\delta\), we obtain that
		\begin{equation}\label{eq:multidigit-square-sum-bound}
			\sum_{k=1}^{\infty}\|g(D_k)y_k u_j\|^{2} \le C_\delta:=-\frac{1}{c}\log \delta.
		\end{equation}
		
		On the other hand, if $\tau$ is an  accumulation point of \(\{y_{p_j}u_j\}\),  then $1\le |\tau|\le B$; after passing to a subsequence if necessary, we may assume that
		\begin{equation*}
			y_{p_j}u_j\to \tau \qquad \text{ as } j\to\infty.
		\end{equation*}
		Furthermore, since  \(g(D_{k})\in \{1,\ldots,B\}\) and \(m_{k}\in \{2,\ldots,B\}\) for all $k$, by Cantor's diagonal argument, there is a subsequence of $p_j$ (for simplicity, we still denoted it by $p_j$)  such that for each integer $r\geq 0$, the two sequences  \(\{g(D_{p_j-r})\}\) and \(\{m_{p_j-r+1}\}\) are eventually constant as $j$ increases. More precisely, to each \(r\ge0\), there exists an integer $j(r)$ such that \(j\ge j(r)\) implies that \[g(D_{p_j-r})=h_r \quad \text{and}\quad m_{p_j-r+1}=n_r\] for some $h_r\in\{1,\ldots,B\}$ and $n_r\in\{2,\ldots,B\}.$
		
		Put \(N_0=1\) and \(N_r=n_1\cdots n_r\) for \(r\ge1\).  Thus, for every fixed integer \(r\ge0\), we have $y_{p_j-r}=m_{p_j-r+1}\cdots m_{p_j}y_{p_j}=N_ry_{p_j},$ which gives
		\[g(D_{p_j-r})y_{p_j-r}u_j\longrightarrow h_rN_r\tau \qquad \text{ as } j\to\infty.\]
		
		Fix a positive integer $M$ and restrict the left-hand  side of \eqref{eq:multidigit-square-sum-bound} to the indices \(k=p_j-r\), where \(0\le r\le M\), we obtain that $\sum_{r=0}^{M}\|g(D_{p_j-r})y_{p_j-r}u_j\|^{2} \le C_\delta$. Letting \(j\to\infty\), we have
		\[\sum_{r=0}^{M}\|h_rN_r\tau\|^{2}\le C_\delta.\]
		Letting \(M\to\infty\), we obtain
		\[\sum_{r=0}^{\infty}\|h_rN_r\tau\|^{2}<\infty,\] which implies  $\|h_rN_r\tau\|\to0$ as $r\to\infty$.
		
		We write
		\[h_rN_r\tau=v_r+\varepsilon_r, \]
		where \(v_r\in\mathbb Z\) and \(|\varepsilon_r|\leq \frac{1}{2}\). Clearly, $|\varepsilon_r|=\|h_rN_r\tau\|\to0$ as $r\to\infty$. A simple calculation shows that \[h_rv_{r+1}-h_{r+1}n_{r+1}v_r=h_{r+1}n_{r+1}\varepsilon_r-h_r\varepsilon_{r+1}\to0 \qquad \text{as }r\to\infty.\]
		Since all $h_r,v_r,n_r$ are integers, the left-hand side is an integer sequence converging to
		$0$, hence it is eventually zero.  Therefore, for all sufficiently large \(r\), we have   \[\frac{v_{r+1}}{h_{r+1}N_{r+1}}=\frac{v_r}{h_rN_r},\]  so the sequence $\{\frac{v_r}{h_rN_r}\}$ is eventually a rational constant. Consequently, \[\tau=\lim_{r\to\infty}\left(\frac{v_r}{h_rN_r}+\frac{\varepsilon_r}{h_rN_r}\right)\in \mathbb Q,\]
		as desired. This finishes the proof.
	\end{proof}
	Now we show that the condition \eqref{eq:product-cross-decay} of Lemma \ref{lem:product-kaufman} holds for the infinite convolution measure $\mu$, outside the exceptional set $\Lambda:=\mathbb Q$.
	\begin{lemma}\label{cor:multidigit-two-sequence}
If the non-zero real numbers  \(s,t\) satisfy that \(\frac{s}{t}\notin\mathbb Q\), then 		\[\widehat\mu(su)\,\widehat\mu(tu)\to0 \qquad (|u|\to\infty).\]
	\end{lemma}
	
	\begin{proof}
		Suppose on the contrary that there exist $\delta>0$ and a sequence $\{u_j\}$ with \(|u_j|\to\infty\) such that \[|\widehat\mu(su_j)|\ge\delta, \qquad |\widehat\mu(tu_j)|\ge\delta\]  for all sufficiently large $j$. By Lemma \ref{lem:multidigit-integer-rigidity}, for each such $j$,  there are indices \(p_j,q_j\) such that
		\begin{equation}\label{eq.bound}
			1\le |y_{p_j}su_j|\leq B, \qquad 	1\le |y_{q_j}tu_j|\leq B.
		\end{equation}
		Now it follows from \(\frac{y_{p_j}}{y_{q_j}}=\frac{y_{p_j}su_j}{y_{q_j}tu_j}\,\frac{t}{s}\)
		that  \[\frac{|t|}{B|s|}\le	\frac{y_{p_j}}{y_{q_j}}\le \frac{B|t|}{|s|}.\]
		Recall that \(y_{k}=m_{k+1}y_{k+1}\)  for all $k$. Thus,
		\[\frac{y_{p_j}}{y_{q_j}}=
		\begin{cases}
			m_{p_j+1}\cdots m_{q_j}, & p_j<q_j,\\[2mm]
			1, & p_j=q_j,\\[2mm]
			\bigl(m_{q_j+1}\cdots m_{p_j}\bigr)^{-1}, & p_j>q_j.
		\end{cases}\]
		It follows from \(2\le m_{k}\le B\) that there is a positive integer \(C_{B,s,t}\), depending only on $B,s$ and $t$,  such that \[|p_j-q_j|\le C_{B,s,t}.\]
		
		Therefore, by extracting a subsequence of $\{j\}$ if necessary, we may suppose that the sequence $\left\{\frac{y_{p_j}}{y_{q_j}}\right\}$ is eventually a constant with the rational number $\rho\in\mathbb Q$, i.e., for $j$ large enough, we have
		\begin{equation}\label{eq:two-sequence-rho}
			\frac{y_{p_j}}{y_{q_j}}=\rho\in\mathbb Q.
		\end{equation}
		
		Based on \eqref{eq.bound}, by passing to a further subsequence of $\{j\}$,  we may assume that
		\[y_{p_j}su_j\longrightarrow\tau_1, \qquad 	y_{q_j}tu_j\longrightarrow\tau_2 \qquad (j\to\infty). \]
		Applying Lemma \ref{lem:multidigit-integer-rigidity} to the two sequences \(\{su_j\}\) and \(\{tu_j\}\) respectively, we get \(\tau_1,\tau_2\in\mathbb Q\).    Together with \eqref{eq:two-sequence-rho}, we obtain
		\[\frac{s}{t}=\frac{1}{\rho}\frac{\tau_1}{\tau_2}\in\mathbb Q,\] contradicting the assumption that $\frac{s}{t}$ is irrational. This finishes the proof.
	\end{proof}

	\begin{proof}[Proof of Theorem \ref{thm:intro-bounded-integer-moran}]
		Suppose that $F$ is uncountable. By Lemma \ref{cor:multidigit-two-sequence}, the natural measure \(\mu\) supported on \(K\) satisfies \[\widehat\mu(u)\,\widehat\mu(\tau u)\to0 \qquad (|u|\to\infty)\]
		for every \(\tau\notin\mathbb Q\).  Thus, by Lemma \ref{lem:product-kaufman}, with exceptional set \(\Lambda=\mathbb Q\), we get that \(F\cdot K\) is an \(M_0\)-set.
		
		Conversely, assume that \(K\) is a \(U\)-set. Thus, if \(F\) is countable, by  Bari's theorem (\!\!\cite[p. 10]{Bari1949}), the set  \(F\cdot K=\bigcup_{a\in F}aK\) is a \(U\)-set, not an $M_0$-set, so the necessity follows. The proof is complete.
	\end{proof}
	
	\begin{remark}
 We remark here that the equivalence in Theorem \ref{thm:intro-bounded-integer-moran} relies heavily on the assumption \(2\le m_k\le B<\infty\). It not only gives uniform estimates for the mask polynomials in Lemma \ref{lem:uniform-mask-estimate}, but also provides the two-scale Fourier estimate for the convolution measure as in  Lemma \ref{cor:multidigit-two-sequence}. On the other hand, the conclusions of  Theorem \ref{thm:unbounded-m0-product} and Theorem \ref{thm:sharpness}  also show that this assumption is indispensable.
\end{remark}
Under the assumption that \(2\le m_k\le B<\infty\), by  Lemma \ref{lem:moran-dimension}\textup{(i)}, the lower bound of the Hausdorff dimension of homogeneous Moran set $K$ satisfies that
\[\dim_H K=\liminf_{n\to\infty}
\frac{\log\bigl(\prod_{k=1}^n\#D_k\bigr)}{\log Q_n}\ge\frac{\log 2}{\log B}>0.\]

\begin{proof}[Proof of Theorem  \ref{rem:bounded-all-dimensions}]
Fix $s\in(0,1]$ and choose an integer $M\ge3$ such that 		$a:=\log 2/\log M\le s$. Let $K$ be the homogeneous Moran set  \eqref{1.12} satisfies the following conditions. Put $m_n=M$ for all $n$, and choose $D_n$ to be either $\{0,1\}$ or $\{0,1,\ldots,M-1\}$ such that the latter choice having asymptotic frequency $\theta=(s-a)/(1-a)$. Thus, $Q_n=M^n$ and $$A_n=\prod_{k=1}^nN_k= M^{\#\{1\leq k\leq n: D_k=\{0,1,\ldots,M-1\}\}} 2^{\#\{1\leq k\leq n: D_k=\{0,1\}\}} $$ By the choice of the digit sets, we have \[\log A_n=n\bigl[\theta\log M+(1-\theta)\log2\bigr]+o(n).\]By Lemma \ref{lem:moran-dimension}\textup{(i)}, we have $$\dim_H K = \lim_{n\to\infty} \frac{\log A_n}{n\log M} = (1 - \theta)\frac{\log 2}{\log M} + \theta=(1-\theta)a+\theta=s.$$

If \(\theta<1\), the choice \(D_n=\{0,1\}\) occurs infinitely often, we let \(\{n_j\}_{j=1}^\infty\) be the increasing sequence of these indices. If \(\theta=1\), we choose these indices to be \(n_j=2^j\) (\(j\ge 1\)). Since the set \(\{2^j\}\) has asymptotic density zero, this choice preserves the prescribed asymptotic frequency.

Notice that any $x \in K$ can be expressed as  the following expansion
\[x = \sum_{k=1}^{\infty} d_k M^{-k}, \quad \text{where } d_k \in D_k \subseteq \{0, 1, \dots, M-1\},\]which gives
\[M^{n_j-1}x = \underbrace{\sum_{k=1}^{n_j-1} d_k M^{n_j-1-k}}_{\text{integer}} + d_{n_j}M^{-1} + \sum_{k=n_j+1}^{\infty} d_k M^{n_j-1-k}.\]Since  \(D_{n_j}=\{0,1\}\), we have \(d_{n_j} \le 1\), and for all \(k \ge n_j+1\), \(d_k \le M-1\). Therefore, for every \(x\in K\),
\[\{M^{n_j-1}x\} = d_{n_j}M^{-1} + \sum_{k=n_j+1}^{\infty} d_k M^{n_j-1-k}\le \frac1M+(M-1)\sum_{\ell=2}^{\infty}M^{-\ell} = \frac2M < 1,\] so $\{M^{n_j-1}x\}(\!\!\!\!\mod 1) \cap (\frac2M,1)=\emptyset$ for all $x\in K.$ Consequently, \(K\) is an \(H\)-set, which immediately implies that \(K\) is a \(U\)-set. Since \(K\) is an uncountable compact set, applying Theorem  \ref{thm:intro-bounded-integer-moran} yields that \(K\cdot K\) is an \(M_0\)-set.
\end{proof}

	\section{On unbounded integer-based homogeneous Moran sets with consecutive digit sets} \label{sec:unbounded-dense}
This section is devoted to proving Theorem \ref{thm:unbounded-m0-product}. In this section, we always assume that \(K\) is the homogeneous Moran set \(K\) defined as in \eqref{1.12} such that  $m_n,N_n,Q_n$ satisfy the conditions of Theorem \ref{thm:unbounded-m0-product}, and assume that $\mu$ is the natural probability measure supported on $K$:  \[\mu=\overset{\infty}{\underset{n=1}{\scalebox{2}{\(\ast\)}}} \left(\frac1{N_n}\sum_{d_n=0}^{N_n-1}\delta_{d_n/Q_n}\right).\]
 	
The following lemma gives a characterization on  the frequencies at which $|\widehat\mu|$ can stay bounded away from zero.
	\begin{lemma}\label{lem:unbounded-peak-localization}
		Assume that there exists $0<c<1$ such that $N_n\ge cm_n$ for every $n$. If $\delta>0$ and $|u_j|\to\infty$ with 	$|\widehat\mu(u_j)|\ge\delta$, then, after passing to a subsequence, there are indices $k_j\to\infty$ and a non-zero integer $a$ such that $\lim_{j\to\infty}\frac{u_j}{Q_{k_j}}= a.$	\end{lemma}
	
	\begin{proof}
		First, since $|\widehat\mu(-u)|=|\widehat\mu(u)|$ for all $u\in\R$, we may without loss of generality assume that $u_j\to+\infty$ such that $|\widehat\mu(u_j)|\ge\delta$,

		Let $D_n=\{0,1,\ldots,N_n-1\}$ and $P_{D_n}(u)=\frac1{N_n}\sum_{d=0}^{N_n-1}e^{-2\pi i du}$. By simple calculation, the Fourier transform of the measure $\mu$ satisfies that
		\[\widehat\mu(u)=\prod_{n=1}^{\infty}P_{D_n}(u/Q_n).\]
		Notice that $P_{D_n}(u)=1$ for all $u\in\mathbb Z$, and the geometric sum formula gives
		\[P_{D_n}(u)=\frac{1-e^{-2\pi iN_nu}}{N_n(1-e^{-2\pi iu})} \qquad \text{whenever } u\notin\mathbb Z,\] 	
		which yields, from $|\sin(\pi u)|\ge2\|u\|$, that
		\[|P_{D_n}(u)|=\frac{|\sin(\pi N_nu)|}{N_n|\sin(\pi u)|}		\le \frac1{2N_n\|u\|}, 	\qquad \|u\|:=\dist(u,\mathbb Z).\]
		Thus, for every $j$ and every $n\ge1$, the assumption $|\widehat\mu(u_j)|\ge\delta$ and $|\widehat\mu(u)|\le |P_{D_n}(u/Q_n)| ~(\forall u\in\R)$ implies
		\begin{equation}\label{eq:unbounded-near-integer}
			\left\|\frac{u_j}{Q_n}\right\|\le\frac1{2\delta N_n}\qquad(n\ge1).
		\end{equation}
		Since $u_j\to+\infty$ and $m_j\to+\infty$, for each $j$ there is a unique integer $r_j\ge1$ such that $Q_{r_j}\le u_j<Q_{r_j+1}$. Thus,
		\begin{equation}\label{in.eq.1}
			\frac1{m_{r_j+1}}\le x_j:=\frac{u_j}{Q_{r_j+1}}<1.
		\end{equation}
		It together with \eqref{eq:unbounded-near-integer} gives
		\begin{equation}\label{in.eq.2}
\min\{x_j,1-x_j\}=\|x_j\|\le \frac1{2\delta N_{r_j+1}}			\le \frac1{2\delta c\,m_{r_j+1}}=: \varepsilon_j\to0 \qquad \text{as } j\to\infty.
		\end{equation}
		By passing to a subsequence of $\{x_j\}_j$ if necessary, we still denote it by $\{x_j\}_j$ and assume that either $x_j\le\varepsilon_j$ or $1-x_j\le\varepsilon_j$ for all $j$.
		
If $1-x_j\le\varepsilon_j$ for every $j$, by taking $k_j=r_j+1$ in \eqref{in.eq.1}, we obtain that	$\frac{u_j}{Q_{k_j}}= x_j\to 1$ as $j\to\infty$, as desired.
		
		If  $x_j\le\varepsilon_j$ for every $j$, we put
		$z_j:=\frac{u_j}{Q_{k_j}}$ where $k_j=r_j$. It follows from \eqref{eq:unbounded-near-integer} that the sequence $\{z_j\}$  satisfies that
		\[\|z_j\|=\left\|\frac{u_j}{Q_{r_j}}\right\|
		\le\frac1{2\delta N_{r_j}}\leq \frac{1}{2\delta c m_{r_j}}\to0, \qquad \text{as } \quad j\to\infty.\]
		On the other hand, by \eqref{in.eq.1} and \eqref{in.eq.2}, we get that  $z_j:=\frac{u_j}{Q_{k_j}}=m_{r_j+1}x_j$ satisfies that 	$1\le z_j\le \frac1{2\delta c}$.
		Now we choose $a_j\in\mathbb Z$ such that $|z_j-a_j|=\|z_j\|$. The sequence $\{z_j\}$ lies in a fixed compact interval, so $\{a_j\}$ takes only finitely many values. After passing to a 		further subsequence, we may take $a_j=a_0$ for all $j$, where $a_0\ge1$, and hence 	$z_j\to a_0$. Finally,  \[\frac{u_j}{Q_{k_j}}=z_j\to  a_0 \in\mathbb Z\setminus\{0\}.\]
		Since $k_j=r_j\to\infty$, this finishes the proof.
	\end{proof}
	
	\begin{lemma}\label{lem:unbounded-ratio}
		Assume that there exists $0<c<1$ such that $N_n\ge cm_n$
		for every $n$. Suppose that $|u_j|,|v_j|\to\infty$,
		$|\widehat\mu(u_j)|,|\widehat\mu(v_j)|\ge\delta>0$, and
		$0<c_0\le|u_j/v_j|\le C_0$.  Then every accumulation point of $u_j/v_j$ is rational.
	\end{lemma}
	
	\begin{proof}
		Let $\lambda$ be an accumulation point of $\{u_j/v_j\}$. Passing to a subsequence if necessary, we may assume that the sequence $u_j/v_j$ converges to $\lambda$. By the assumption, we have  $c_0\le |\lambda| \le C_0$. 	According to  Lemma \ref{lem:unbounded-peak-localization},  there exist  		sequences $k_j,\ell_j\to \infty$ and
		integers $a,b\in\mathbb Z\setminus\{0\}$ such that 		$u_j/Q_{k_j}\to a$ and $v_j/Q_{\ell_j}\to b$.
		Therefore,
		\[\frac{Q_{k_j}}{Q_{\ell_j}}=\frac{u_j}{v_j}
		\frac{v_j/Q_{\ell_j}}{u_j/Q_{k_j}}
		\to \lambda\frac ba, \qquad \text{as }~ j\to\infty,\] where the limit is finite and nonzero. However, it is impossible, since $Q_n=m_1m_2\ldots m_n$, we have
		\begin{itemize}
			\item  if $k_j>\ell_j$ along an  infinite subsequence, then 		\[\frac{Q_{k_j}}{Q_{\ell_j}}=m_{\ell_j+1}\cdots m_{k_j}
			\ge m_{\ell_j+1}\to\infty,\]
			\item if $k_j<\ell_j$ along an infinite subsequence, then 		$Q_{k_j}/Q_{\ell_j}\le m_{k_j+1}^{-1}\to0$.
		\end{itemize}
		Therefore, $k_j=\ell_j$ for all sufficiently large $j$, and hence 	\[\lambda=\lim_j\frac{u_j}{v_j}=\lim_j\frac{u_j/Q_{k_j}}{v_j/Q_{k_j}}
		=\frac ab\in\mathbb Q.\]
		This finishes the proof.
	\end{proof}

	The following lemma gives a simple  criterion for the uniqueness of homogeneous Moran sets with consecutive digit sets, which will  be used in the proof of Theorem \ref{thm:unbounded-m0-product}.
	\begin{lemma}\label{lem:consecutive-digit-h}
		Let $m_k\ge2$ and $1\le N_k\le m_k$ be integers with $N_k\le Lm_k$ for every $k\ge1$, where $0<L<1$. Put $Q_0=1$, $Q_n=m_1\cdots m_n$, and consider the homogeneous Moran set \[K=\left\{\sum_{k\ge1}\frac{d_k}{Q_k}:d_k\in\{0,1,\ldots,N_k-1\}\right\}.\]
		Then $K\subset[0,L]$ is a $H$-set, and hence a $U$-set.
		Moreover,  \begin{equation}\label{eq:consecutive-tail-upper}
R_n:=\sum_{k>n}\frac{N_k-1}{Q_k}\le\frac{N_{n+1}}{Q_{n+1}}			\le\frac{L}{Q_n}\qquad(n\ge0).
		\end{equation}
	\end{lemma}
	\begin{proof}
Since $(m_k-1)/Q_k=Q_{k-1}^{-1}-Q_k^{-1}$ and  $N_k\le Lm_k$ for every $k\ge1$, it follows that
\[\begin{aligned}
R_n&\le\frac{N_{n+1}-1}{Q_{n+1}}+\sum_{k>n+1}\frac{m_k-1}{Q_k} =\frac{N_{n+1}-1}{Q_{n+1}}+\sum_{k>n+1}\left(\frac1{Q_{k-1}}-\frac1{Q_k}\right)
=\frac{N_{n+1}-1}{Q_{n+1}}+\frac1{Q_{n+1}}			=\frac{N_{n+1}}{Q_{n+1}}\le\frac{L}{Q_n}.
\end{aligned}\]
In particular, taking $n=0$ yields that $K\subset[0,L]$.

Therefore, for $x\in K$, it is not hard to check that  $\{Q_nx\}\notin(L,1)$ for all $x\in K$ and $n\ge1$. Since $Q_n$ is a strictly increasing sequence of positive integers, $K$ is an $H$-set. Therefore, by Theorem \ref{thm:salem-zygmund-criterion}, $K$ is a $U$-set. The proof is finished.
	\end{proof}
	
	We are ready to prove Theorem \ref{thm:unbounded-m0-product}.
	\begin{proof}[Proof of Theorem \ref{thm:unbounded-m0-product}] 	
	First, we prove the statement (a).	Assume that $\inf_n N_n/m_n=0$. Let $\varepsilon_n=N_n/m_n$ and $R_n=\sum_{k>n}\frac{N_k-1}{Q_k}.$

		We construct the Cantor set $F$ as follows.  Let $\eta_j=2^{-j-4}$,  $p_0=n_0=1$. We start with 	$\mathcal I_0=\{[1/4,1/2]\}$, and construct, by recursion, finite families $\mathcal I_j$ of $2^j$  pairwise disjoint non-degenerate closed intervals.  Given $\mathcal I_{j-1}$, where $j\geq1$, we choose an integer $p_j>p_{j-1}$ sufficiently large such that every $I \in \mathcal I_{j-1}$ contains  two distinct points $\frac{a_{j,1}}{p_j}, \frac{a_{j,2}}{p_j},$ where $a_{j,1}\not= a_{j,2}\in\mathbb Z.$
		Moreover, since $\inf_n\varepsilon_n=0$, we can choose $n_j>n_{j-1}$ such that 	 \[\varepsilon_{n_j}\le\frac{\eta_j}{2p_j},\qquad \text{and}\qquad  q_j:=p_jQ_{n_j-1}.\]
		Finally, choose a common length $\ell_j>0$ sufficiently
		small that the intervals
		\[\left[\frac {a_{j,1}}{p_j},\frac {a_{j,1}}{p_j}+\ell_j\right]\subset I\qquad \text{and}\qquad
		\left[\frac {a_{j,2}}{p_j},\frac {a_{j,2}}{p_j}+\ell_j\right]\subset I,\] where
		\[\ell_j\le \min\left\{\frac{\eta_j}{2q_j},~ 2^{-j^2}\right\}.\]
		Define the Cantor set as
		\[F=\bigcap_{j\ge1}\bigcup_{I\in\mathcal I_j}I.\]
		Thus $F\subset[1/4,1/2]$ is a nonempty Cantor set, which is uncountable. Note that, for every $t>0$,
		\[\sum_{I\in\mathcal I_j}|I|^t\le2^{j-tj^2}\longrightarrow0 \qquad (j\to\infty).\]
		Consequently, $\dim_HF=0$.
		
 Notice that for every $j\ge1$, each $x\in K $ admits the expansion $x=\sum_{k=1}^{n_j-1}\frac{d_k}{Q_k}+\sum_{k=n_j}^{\infty}\frac{d_k}{Q_k}$, where $d_k\in\{0,1,\ldots,N_k-1\}$.   By \eqref{eq:consecutive-tail-upper}, we can write $Q_{n_j-1}x$ as
$$Q_{n_j-1}x=a_{n_j}(x)+\theta_{n_j}(x),\qquad a_{n_j}(x)\in\mathbb Z,\qquad
0\le\theta_{n_j}(x)\le Q_{n-1}R_{n-1}	\le\varepsilon_{n_j}.$$
 Without loss of generality, we can write each $y\in F$ as  $$y=\frac b{p_j}+r,\qquad \qquad b\in\mathbb N,~ 0\le r\le\ell_j.$$
Therefore, for the strictly increasing sequence $\{q_j\}_j$ of positive integers, we have 		\[q_jxy=ba_{n_j}(x)+b\theta_{n_j}(x)+q_jxr.\]
Since $0\le x\le1$ and $0\le b\le p_j$, we obtain 	\[0\le q_jxy-ba_{n_j}(x)\le p_j\varepsilon_{n_j}+q_j\ell_j 	\le\eta_j<\frac1{20}.\]
		As $ba_{n_j}(x)\in\mathbb Z$, it follows immediately that
\[\{q_jz\}\notin(1/20,1)\qquad(z\in F\cdot K ,\ j\ge1).\]
		Therefore the compact set $F\cdot K \subset[0,1]$ is an $H$-set. By Theorem \ref{thm:salem-zygmund-criterion}, it is a $U$-set, and hence is not an $M_0$-set. The statement (a) is proved.

\smallskip

	Second, we prove the statement (b).
 Fix $s \in (\frac12, 1]$. Let \(K\) be the homogeneous Moran set defined as in \eqref{1.12} such that \(m_n=2^{n^2+n}\) and the digit sets \(D_n = \{0, \ldots, N_n-1\}\), where \(N_n\) satisfies the following conditions:
\begin{itemize}
    \item If \(s=1\), we take \(N_n = 2^{n^2}\). In this case, \(N_n/m_n \to 0\) and  $\frac{N_n}{m_n}=2^{-n}\leq \frac12$ for all $n$.
    \item If \(s \in (1/2, 1)\), we take \(N_n = \lfloor 2^{s(n^2+n)} \rfloor\). In this case, \(N_n/m_n \to 0\) and   $\frac{N_n}{m_n}\leq2^{-(1-s)(n^2+n)} \leq 4^{-(1-s)}<1$ for all $n$.
\end{itemize}
Therefore,  \(K\) is a \(U\)-set by Lemma \ref{lem:consecutive-digit-h}.

We next compute the Hausdorff dimension of \(K\).  If $s=1$, since
\[\begin{aligned}
\log A_n&=\Bigl(\sum_{k=1}^n k^2\Bigr)\log 2, \qquad
-\log\frac{N_{n+1}}{Q_{n+1}} =\log Q_n+\log\frac{m_{n+1}}{N_{n+1}}
=\Bigl(\sum_{k=1}^n(k^2+k)+n+1\Bigr)\log 2,
\end{aligned}\]
it follows from  Lemma \ref{lem:moran-dimension}\textup{(ii)} that \[\dim_H K= \lim_{n\to\infty} \frac{\log A_n}{-\log\frac{N_{n+1}}{Q_{n+1}}}=\liminf_{n\to\infty}\frac{\sum_{k=1}^n k^2}{\sum_{k=1}^n(k^2+k)+n+1}=1.\]

If $s \in (\frac12, 1)$, it follows from  \(N_k = \lfloor 2^{s(k^2+k)} \rfloor\)  that \(\log N_k = s(k^2+k)\log 2 + O(1)\). It together with the expression \(\sum_{k=1}^n (k^2+k) = \frac{n(n+1)(n+2)}{3}\) gives
\[\log A_n = \sum_{k=1}^n \log N_k = s \sum_{k=1}^n (k^2+k)\log 2 + O(n)= s \frac{n(n+1)(n+2)}{3} \log 2 + O(n),\]
and
\[\begin{aligned}
-\log\frac{N_{n+1}}{Q_{n+1}}
&=\log Q_{n+1} - \log N_{n+1}
=\sum_{k=1}^{n+1} (k^2+k)\log 2 - \left( s(n+1)(n+2)\log 2 + O(1) \right) \\
&= \frac{n(n+1)(n+2)}{3} \log 2 + (1-s)(n+1)(n+2)\log 2 + O(1).
\end{aligned}
\]
Therefore, by Lemma \ref{lem:moran-dimension}\textup{(ii)},
\[\begin{aligned}
\dim_H K &= \lim_{n\to\infty} \frac{\log A_n}{-\log\frac{N_{n+1}}{Q_{n+1}}}=s.
\end{aligned}\]
In either case, by Theorem \ref{thm.OSS}, $K \cdot K $ is an $M_0$-set. Thus, the statement (b) is proved.

 		Third, we show the equivalence of (i) and (ii).

(ii) $\Rightarrow$ (i). It follows from the statement (a).

 (i) $\Rightarrow$ (ii). 
 If (i) holds, then there is a constant $0<c<1$ such that $N_n\ge cm_n$ for every $n$. Let $\mu$ be its natural probability measure.
We claim that $$\text{$\lim_{|u|\to\infty}\widehat\mu(u)\widehat\mu(\tau u)=0$ for all $\tau\notin\mathbb Q$.}$$

 In fact, if there exists a $\delta>0$ and a sequence $|u_j|\to\infty$ such that $|\widehat\mu(u_j)|,|\widehat\mu(\tau u_j)|\ge\delta$, by applying Lemma \ref{lem:unbounded-ratio} to the two sequences $\{u_j\}$ and $\{\tau u_j\}$, we obtain that $1/\tau$ is rational, which is a contradiction. The claim is proved.

 Therefore,  Lemma \ref{lem:product-kaufman}, with $\Lambda:=\mathbb Q$, yields that $F\cdot K $ is an $M_0$-set for every uncountable compact set $F\subset\mathbb R$, (ii) follows.

\smallskip
If $N_n\le Lm_n$ for every $n$, where $0<L<1$, by Lemma \ref{lem:consecutive-digit-h}, $K $ is an $H$-set, and hence a $U$-set.

For the proof of \eqref{eq.equiv}, the sufficiency follows from the statement (ii); conversely, if $F$ is countable, then $F\cdot K=\bigcup_{a\in F}aK$ is a countable union of closed $U$-sets, by the Marcinkiewicz--Zygmund theorem, with $a=0$ treated separately. By Bari's theorem, $F\cdot K$ is a $U$-set, this is a contradiction; the necessity follows.

Finally, we prove that $\dim_HK =1$. Indeed, since  $A_n=N_1N_2\ldots N_n, Q_n=m_1m_2\ldots m_n$ with $m_n\to\infty$, it follows from $cm_k\le N_k\le m_k$ that  $c^nQ_n\le A_n\le Q_n$, which is equivalent to saying that $$n\log c +\log Q_n\le \log A_n\le \log Q_n.$$  It together with $m_n\to\infty$ yields that $n/\log Q_n\to0$ and hence $\log A_n/\log Q_n\to1$. On the other hand, since $N_{n+1}\ge cm_{n+1}$, for all sufficiently large $n$ we  have 	\[R_n\ge\frac{N_{n+1}-1}{Q_{n+1}} \ge\frac{c}{2Q_n}.\] Combining this with $R_n\le\frac{N_{n+1}}{Q_{n+1}}\le Q_n^{-1}$ as in \eqref{eq:consecutive-tail-upper}, we get that
$$\lim_{n\to\infty} \frac{\log R_n}{\log \frac{N_{n+1}}{Q_{n+1}}}=\lim_{n\to\infty}\frac{\log R_n}{-\log Q_n}=1.$$
By Lemma \ref{lem:moran-dimension}\textup{(ii)}, $$\dim_HK =\liminf_{n\to\infty}\frac{\log A_n}{-\log \frac{N_{n+1}}{Q_{n+1}}}=1.$$
This finishes the proof of Theorem \ref{thm:unbounded-m0-product}.
	\end{proof}
	
Based on Theorems \ref{thm:intro-bounded-integer-moran} and \ref{thm:unbounded-m0-product}, we close this section with the following question.
	\begin{question}\label{que:unbounded-moran}
		For a homogeneous  Moran set $K$ as in \eqref{1.12} generated by a sequence of unbounded bases \(\{m_k\}_{k\ge1}\) and digit sets \((D_k)_{k\ge1}\), what additional assumptions ensure that the natural measure $\mu$ as in \eqref{eq.4.1} satisfies
\[\widehat\mu(u)\widehat\mu(\tau u)\to0\qquad (|u|\to\infty)\]		for every \(\tau\) outside a countable exceptional set? Furthermore, under what conditions does the following conclusion
$$F\cdot K \text{ is an }M_0\text{-set}		\quad\text{if and only if}\quad F\text{ is uncountable}$$
remain valid when the bases \(\{m_k\}_{k\ge1}\) are unbounded?
	\end{question}
	
	\section{On sparse integer-based homogeneous Moran sets and the $1/2$ bound}\label{sec:sharpness}
	
	The goal of this section is to prove Theorem \ref{thm:sharpness}.
	
	\subsection{Proof of Theorem \ref{thm:sharpness}\textup{(a)}}
For the proof of Theorem \ref{thm:sharpness}(a), we require a construction of homogeneous Moran sets which, unlike the approach in \cite[Theorem 2.2]{FengWenWu1997}, also establishes a new {\it dimension intermediate value theorem} for such sets.

	\begin{lemma}\label{lem:two-scale-frostman}
Fix a real number $\alpha\ge2$ and a sequence of integers  $\{b_j\}_{j=1}^\infty$ such that $b_j\ge\lceil\alpha\rceil$ for all $j$ and $b_j\to\infty$. Let $m_1\geq 400$ and  $N_1=\left\lfloor\frac{m_1}{100}\right\rfloor.$ For $k\geq1$, put $$Q_k=m_1\cdots m_k,\quad m_{k+1}=Q_k^{b_k}, \quad  \text{and}\quad  N_{k+1}=\left\lfloor\frac{m_{k+1}}{100 Q_k^{\alpha-1}}\right\rfloor.$$
 Consider the following homogeneous Moran set as in \eqref{1.12}:
 \[K :=\left\{\sum_{j=1}^{\infty}\frac{d_j}{Q_j}: 	d_j\in\{0,1,\ldots,N_j-1\},~j\geq1\right\}.\]
Let $A_j=N_1\cdots N_j$  and $R_j=\sum_{k>j}\frac{N_k-1}{Q_k}$ for $j\ge0$. Then
\begin{equation}\label{eq:two-scale-data}
\frac1{200Q_j^\alpha}\le R_j\le\frac1{75Q_j^\alpha} \quad(j\ge1),
			\qquad 	\lim_{j\to\infty}\frac{\log A_j}{\log Q_j}=1.
		\end{equation}
The natural probability measure
$\mu:=\overset{\infty}{\underset{n=1}{\scalebox{2}{\(\ast\)}}}
\left(\frac1{N_n}\sum_{d_n=0}^{N_n-1}\delta_{d_n/Q_n}\right)$, supported on the set $K$, satisfies the Frostman condition. More precisely, for every \(0<t<1/\alpha\), there exist constants \(C>0\) and \(\varepsilon_t>0\) such that \begin{equation}\label{eq:frostman-estimate}
\mu(I)\le C|I|^t
\end{equation}
whenever $I\subset\mathbb R$ is a set with $0<|I|\leq \varepsilon_t$. Consequently, $$\dim_HK=\frac{1}{\alpha}.$$
	\end{lemma}
	
	\begin{proof}
Obviously, $4\leq N_k<m_k$ for all $k\geq1,$ which yields, from $\alpha\geq2$ and the definitions of $N_k,m_k,Q_k$, that
$$\frac{Q_{k+1}}{200 Q_k^{\alpha}}=\frac{m_{k+1}}{200 Q_k^{\alpha-1}} \leq\frac{m_{k+1}}{100 Q_k^{\alpha-1}}-2<N_{k+1}-1<\frac{m_{k+1}}{100 Q_k^{\alpha-1}}=\frac{Q_{k+1}}{100 Q_k^{\alpha}}.$$
Combining this with $Q_{k+1}=Q_k^{b_k+1}>2Q_k$ and $\alpha\ge2$, we get that	\[\frac1{200Q_j^\alpha}<\frac{N_{j+1}-1}{Q_{j+1}}	<R_j=\sum_{k>j}\frac{N_k-1}{Q_k}\le\frac1{100Q_j^\alpha}\sum_{r=0}^{\infty}2^{-\alpha r}	\le\frac1{75Q_j^\alpha}.\]
		
		Let $\eta_j=\log A_j/\log Q_j, j\geq1.$ It is not hard to see that $A_j<Q_j$ and $0\le\eta_j<1$. Note that $\frac{m_{k+1}}{Q_k^{\alpha-1}}=Q_{k}^{b_k+1-\alpha}$. With easy calculations,
\[(b_j+1-\alpha)\log Q_j-\log 200 < \log N_{j+1}\leq(b_j+1-\alpha)\log Q_j-\log 100.\] 		Write $\log N_{j+1}=(b_j+1-\alpha)\log Q_j+O(1)$. 		Hence \[\eta_{j+1}=\frac{\log A_j+\log N_{j+1}}{(b_j+1)\log Q_j}=\frac{\eta_j+(b_j+1-\alpha)}{b_j+1}		+O\!\left(\frac1{(b_j+1)\log Q_j}\right),\]
		i.e.,
		$$1-\eta_{j+1}=\frac{\alpha-\eta_j}{b_j+1}		+O\!\left(\frac1{(b_j+1)\log Q_j}\right).$$	Therefore, there is a constant $C$ such that 		\[0<1-\eta_{j+1}\le\frac{\alpha}{b_j+1} 		+\frac{C}{(b_j+1)\log Q_j}\to0 \qquad (j\to\infty).\]
		This proves \eqref{eq:two-scale-data}.
		
For the proof of \eqref{eq:frostman-estimate}, we first note that the set $K$ can be represented as
$$K=\bigcap_{j=1}^\infty \bigcup_{d_i\in \{0,1,\ldots,N_i-1\},1\leq i\leq j} I(d_1,\ldots,d_j),$$
where \[ I(d_1,\ldots,d_j)=\left[\sum_{k\le j}\frac{d_k}{Q_k},\ \sum_{k\le j}\frac{d_k}{Q_k}+R_j\right]\] is the $j$-order fundamental intervals of $K$, whose $\mu$ measure is $\frac{1}{A_j},$ i.e., $$\mu(I(d_1,\ldots,d_j))=\frac{1}{A_j}.$$

		Fix $0<t<1/\alpha$ and choose $\beta$ such that $\alpha t<\beta<1$. By $\lim_{j\to\infty}\frac{\log A_j}{\log Q_j}=1$, there is an integer $j_0$,  depending on $t$, such that $j\geq j_0$ implies that
		\begin{equation}\label{add.eq.1}
			\frac{\log A_j}{\log Q_j}\geq \beta \quad \Leftrightarrow \quad A_j\geq Q_j^\beta.
		\end{equation}
Let  $\varepsilon_t:=R_{j_0}$  and let $I$ be any set with diameter $r$ (i.e., $r:=|I|$) such that $0<|I|\leq\varepsilon_t$. Then there is a unique integer $j\geq j_0+1$ such that
		\begin{equation}\label{add.eq.2}
			R_j< r:=|I|\leq R_{j-1}.
		\end{equation}
Note that all the left endpoints of $j$-order fundamental intervals belong to $Q_j^{-1}\mathbb Z$ and are distinct. If such an interval intersects $I$, then its left endpoint lies in an interval of length at most $r+R_j\le2r$. Hence $I$ intersects at most $1+2rQ_j$ level-$j$ fundamental intervals, giving
		\begin{equation}\label{add.eq.3}
			\mu(I)\leq \frac{1+2rQ_j}{A_j}.
		\end{equation}
		
		\textbf{Case 1.} If $rQ_j\le1$, it follows from $t<\frac{\beta}{\alpha}<\frac{1}{\alpha}\leq \frac12$ that
\[\mu(I)\le\frac3{A_j}\overset{\eqref{add.eq.1}}{\le} 3Q_j^{-\beta}\overset{\eqref{eq:two-scale-data}}{\le} 3(200R_j)^{\beta/\alpha}\overset{\eqref{add.eq.2}}{\le}3(200r)^{\beta/\alpha}\le 3\sqrt{200} \cdot r^t.\]
		
		\textbf{Case 2.} If $rQ_j>1$, it follows from  the inequalities  $N_{j}\ge\frac{Q_{j-1}^{b_{j-1}+1-\alpha}}{200}$,  $A_{j-1}\geq Q_{j-1}^\beta$ (see \eqref{add.eq.1}) and $Q_{j}=Q_{j-1}^{b_{j-1}+1}$ that
		\[\frac{Q_j}{A_j}=\frac{Q_j}{A_{j-1}N_j} \le200Q_{j-1}^{\alpha-\beta}.\] Thus, \eqref{add.eq.3}  together with $r\leq R_{j-1}\le\frac1{75Q_{j-1}^\alpha}$ (by \eqref{add.eq.2} and \ref{eq:two-scale-data}) gives
		\begin{align*}
			\mu(I)  &\leq \frac{3rQ_j}{A_j}\leq 600rQ_{j-1}^{\alpha-\beta}=600r^t \cdot (r^{1-t}Q_{j-1}^{\alpha-\beta})\\
&\le 600r^t\cdot (75Q_{j-1}^\alpha)^{t-1}Q_{j-1}^{\alpha-\beta}
=600 \cdot 75^{t-1}\cdot r^t Q_{j-1}^{\alpha t-\beta}\\
& \leq  600 \cdot 75^{\frac{1}{\alpha}-1}\cdot r^t \qquad (\text{by $\alpha t-\beta<0$ and $ t<\frac{1}{\alpha}$})
		\end{align*}
		This proves  \eqref{eq:frostman-estimate}  by taking $$C=\max\{3\sqrt{200}, 600 \cdot 75^{\frac{1}{\alpha}-1}\}.$$
		
		By the mass distribution principle \cite[Principle 4.2]{Falconer} and  \eqref{eq:frostman-estimate}, we have $\dim_HK\geq t$. Since $0<t<1/\alpha$ is arbitrary, we get  $\dim_HK\ge1/\alpha$.
		Conversely, since the $j$-order fundamental intervals $\{I(d_1,\ldots,d_j): d_i\in \{0,1,\ldots,N_i-1\},i\leq j\}$  forms  a $\delta_j$-covers of $K$ with $\delta_j:=R_j\leq \frac1{75Q_j^\alpha}\to0.$ Recall that $A_j<Q_j$ for all $j$. Thus, for every $s>1/\alpha$, the $s$-dimensional Hausdorff measure $\mathcal H^s$ satisfies that
\[\mathcal H^s(K)\le\liminf_{j\to\infty}A_jR_j^s 		\le75^{-s}\lim_{j\to\infty}Q_j^{1-\alpha s}=0.\]
		Thus $\dim_HK\leq s$. Since $s$ is arbitrary, we have  $\dim_HK\leq 1/\alpha$. Consequently, $\dim_HK=1/\alpha$.
	\end{proof}
	
	\begin{proof}[The proof of Theorem \ref{thm:sharpness}\textup{(a)}]
		Suppose first that $0<s\le1/2$. Put $\alpha=1/s\ge2$ and define a sequence $\{b_j\}_{j=1}^\infty$ such that  $b_j=\lceil\alpha\rceil+j-1$. Let  $E=K$ be the homogeneous Moran set constructed as in Lemma \ref{lem:two-scale-frostman} such that $\dim_HK=1/\alpha.$
Since the integer bases of $K$ satisfy	$m_1=Q_1\geq400$ and $m_{j+1}=Q_j^{b_j}\to\infty$, it follows that
\[\frac{\max D_1}{m_1}\le\frac1{100},\qquad\frac{\max D_{k+1}}{m_{k+1}}\le\frac1{100Q_k^{\alpha-1}}\le\frac1{100}
\quad\text{for all } k.\]
Thus, $E$ is an $H$-set, and hence a $U$-set by Lemma \ref{lem:left-half-moran-u}. It is not hard to see from \eqref{eq:two-scale-data} that
		$E\subset[0,(N_1-1)/Q_1+R_1]\subset[0,1]$.
		
We next prove that $E\cdot E$ is an $H$-set. In fact, for $x,y\in E$, we can without loss of generality write	\[x=\frac{M_j}{Q_j}+r_j,\qquad y=\frac{M'_j}{Q_j}+s_j,\]
		where $M_j,M'_j$ are integers in $[0,Q_j]$ and 	$0\le r_j,s_j\le R_j$. Then \[Q_j^2xy=M_jM'_j+	\bigl(Q_jM_js_j+Q_jM'_jr_j+Q_j^2r_js_j\bigr).\]
Since $\alpha\ge2$ and $R_j\le(75Q_j^2)^{-1}$ (see \eqref{eq:two-scale-data}), we have
\[0\le Q_jM_js_j+Q_jM'_jr_j+Q_j^2r_js_j\le2Q_j^2R_j+Q_j^2R_j^2		\le\frac2{75}+\frac1{75^2Q_j^2}<\frac1{20}. \]
		Therefore,
		\begin{equation*}
\{Q_j^2z\}\notin(1/20,1)\qquad(z\in E\cdot E,\ j\ge1).
		\end{equation*}
This means that $E\cdot E$ is an $H$-set, and Theorem \ref{thm:salem-zygmund-criterion} shows that $E\cdot E$ is a $U$-set.
		
		It remains to consider $s=0$. Let $E_*$ be the $U$-set constructed in Lemma \ref{lem:two-scale-frostman}  under the conditions $\alpha=2$, $b_j=j+1$, $m_1=Q_1$ and $Q_{j+1}=Q_j^{b_j+1}=Q_j^{j+2}$. We consider a subset $E_0$ of $E_*$ as follows:
\[E_0=\left\{\sum_{j\ge1}\frac{d_j}{Q_j}: d_j\in\{0,1\}\right\}\subset E_*.	\]
Note that the $2^j$ fundamental intervals of order $j$ form  a $\delta_j$-covers of $E_0$ with $\delta_j:=R_j,$ which is controlled by a geometric series
\[R_j=\sum_{k>j}\frac1{Q_k} \le \frac{1}{Q_{j+1}} \sum_{n=0}^{\infty} \left( \frac{1}{Q_{j+1}^{j+2}} \right)^n  = \frac{1}{Q_{j+1}} \cdot  \frac{1}{1-Q_{j+1}^{-(j+2)}}\le\frac2{Q_{j+1}}\to 0\qquad (j\to\infty).\]
 Because $Q_{j+1}=Q_j^{j+2}=\cdots=Q_1^{\frac{(j+2)!}{2}}$,  for each $t>0$, the  $t$-dimensional Hausdorff measure of $E_0$  satisfies that \[\mathcal H^t(E_0) \le\liminf_{j\to\infty}2^j(R_j)^t		\le2^t\lim_{j\to\infty}2^jQ_{j+1}^{-t}=0,\]
which gives $\dim_HE_0=0$. As in the proof for the case $0<s\le1/2$, we can similarly show that   $E_0\cdot E_0$  is an $H$-set, and hence a $U$-set. Remember that the class of $U$-sets is closed under taking subsets. Thus, in this case, we can take $E=E_0$.
		
Finally, for any non-zero $x_0\in E\setminus\{0\}$, we have $x_0E\subset E\cdot E$. The scaling invariance of Hausdorff dimension gives \[\dim_H(E\cdot E)\ge\dim_H(x_0E)=\dim_HE=s.\]
This finishes the  proof  of Theorem \ref{thm:sharpness}\textup{(a)}. 	\end{proof}
	
  	\subsection{Proof of Theorem \ref{thm:sharpness}\textup{(b)}}
	We first use the $H$-set condition to obtain a nontrivial Fourier coefficient.	
	\begin{lemma}\label{lem:avoidance-fourier}
		Let \(I\subset(0,1)\) be a nonempty open interval. Then there exist constants \(c_0=c_0(I)>0\) and \(M=M(I)\in\mathbb N\) such that the following holds. For every positive integer \(q\) and every Borel probability measure \(\nu\) on \(\mathbb R\) satisfying $\{qz\}\notin I~ (z\in \operatorname{supp}\nu)$, there exists an integer \(k\) with $0<|k|\le M$ such that $|\widehat{\nu}(kq)|\ge c_0.$
	\end{lemma}
	
	\begin{proof}
		Let $\mathbb T:=\mathbb R/\mathbb Z=[0,1)$ be the circle and let \(\psi\in C^\infty(\mathbb T)\) be a nonnegative smooth function with $\int_0^1\psi(x)\,dx=1$ and $\operatorname{supp} (\psi) \subset I.$ Recall that the Fourier transform of $\psi$ is defined by 	\[\widehat{\psi}(k)=\int_0^1\psi(x)e^{-2\pi ikx}\,dx.\]
		Then the Fourier series $\sum_{k\in\mathbb Z} \widehat{\psi}(k) e^{2\pi ikx}$ of \(\psi\in C^\infty(\mathbb T)\) is absolutely convergent, so
		\[0<C_\psi:=\sum_{k\in\mathbb Z\setminus\{0\}}
		|\widehat{\psi}(k)|<\infty .\]
		This yields that there is an integer \(M\ge1\) sufficiently large so that $\sum_{|k|>M}|\widehat{\psi}(k)|\le\frac12.$
    Now the assumption and the Fourier transform of $\nu$ gives  \[0=\int_{\mathbb R}\psi(qz)\,d\nu(z) =1+\sum_{k\ne0}\widehat{\psi}(k)\widehat{\nu}(-kq).\]
		Since \(\nu\) is a probability measure, it follows that $|\widehat{\nu}(\xi)|\le1$, which gives
		\[1\le \sum_{0<|k|\le M}|\widehat{\psi}(k)||\widehat{\nu}(-kq)|+\sum_{|k|>M}|\widehat{\psi}(k)|.\]
		Therefore, $\sum_{0<|k|\le M}|\widehat{\psi}(k)||\widehat{\nu}(-kq)|\geq \frac12,$ and hence  there exists some integer \(k\) with $0<|k|\le M$
		such that \[|\widehat{\nu}(-kq)| \ge \frac{1}{2\sum_{0<|k|\le M}|\widehat{\psi}(k)|}\ge \frac{1}{2C_\psi}.\]
		Finally, since \(\nu\) is a positive Borel measure,
		$\widehat{\nu}(-\xi)=\overline{\widehat{\nu}(\xi)},$ so we get the desired result by taking $c_0=\frac1{2C_\psi}.$ The proof of the lemma is finished.
	\end{proof}
	
	We also use the following classical large sieve inequality.
	\begin{lemma}[Principle of the large sieve; {\cite[p.~559, Theorem 3]{Montgomery1978}}]
		\label{lem:large-sieve}
If $x_1,\ldots,x_R$ are $\Delta$-separated modulo one, i.e.
		$\|x_r-x_{r'}\|\ge\Delta$ for $r\ne r'$, where 		$\|t\|=\dist(t,\mathbb Z)$, and $a_0,\ldots,a_N\in\mathbb C$, then \[\sum_{r=1}^{R}\Bigl|\sum_{n=0}^{N}a_n \exp(2\pi i n x_r)\Bigr|^2 \ \le\ \bigl(N+\Delta^{-1}\bigr)\sum_{n=0}^{N}|a_n|^2.\]
	\end{lemma}

	The following lemma shows that the existence of sufficiently large Fourier coefficients at prescribed frequencies yields an upper bound for the sum of the lower box dimensions of two compact sets. It is worthy noting that the lemma itself does not require $S$ to be an $H$-set.
	\begin{lemma}\label{lem:fourier-window-product}
		Let $S\subset[0,1]$ be a compact set, $c_0>0$ and $1\le L_j\le U_j$, where $L_j\to\infty$. Assume that for every Borel probability measure $\nu$ supported on $S$ and every $j$, there is a real number $u_j$ such that
		\[L_j\le|u_j|\le U_j, \qquad |\widehat\nu(u_j)|\ge c_0.\]
		If $E_1,E_2\subset[0,1]$ are nonempty compact sets and
		$E_1\cdot E_2\subset S$, then \[\underline\dim_BE_1+\underline\dim_BE_2 \le\liminf_{j\to\infty}\frac{\log U_j}{\log L_j}.\]
	\end{lemma}
	
	\begin{proof}
		Let $E_1,E_2\subset[0,1]$ be nonempty compact sets satisfying $E_1\cdot E_2\subset S$. Set
\begin{equation*}
  \delta_j=\frac{2\varepsilon}{L_j} \qquad \text{with}\qquad 		\varepsilon=\min\left\{\frac18,\frac{c_0}{12\pi}\right\}.
\end{equation*}
		For each $i=1,2$, let $T_i\subset E_i$ be an arbitrary
		nonempty finite $\delta_j$-separated set, i.e., $|t-s|\ge\delta_j$ for all distinct $t,s\in T_i$.
 Let $k_i:=\#T_i$ be the number of the elements in the set $T_i$,  and let \[\nu=\frac1{k_1k_2}\sum_{t\in T_1}\sum_{t'\in T_2}\delta_{tt'}\]
		be the discrete measure supported on the product set $T_1\cdot T_2\subset S$.

 Fix $j$. According to the assumption,  there exists a $u:=u_j\in\mathbb R$ such that $L_j\le |u|\le U_j$ and
		\begin{equation}\label{eq:bilinear-large}
			|\widehat{\nu}(u)|=\left|\frac1{k_1k_2}\sum_{t\in T_1}\sum_{t'\in T_2}e^{-2\pi iutt'}\right|\ge c_0.
		\end{equation}
		
		Let $\kappa:=\frac{\varepsilon}{|u|}$ and $\theta:=u\kappa^2$.  It is clear that  		
\[\kappa\le\frac{\varepsilon}{L_j}=\frac{\delta_j}{2} \qquad\text{and}\qquad  \theta=\frac{\varepsilon^2}{u}.\]
		Given that $T_1$ and $T_2$ are both $\delta_j$-separated sets, we can decompose each $t\in T_1\cup T_2$  as follows:
		\[t=\kappa a_t+r_t, \qquad a_t=\left\lfloor\frac{t}{\kappa}\right\rfloor, \qquad
		0\le r_t<\kappa,\]
		where $a_t\in\{0,\ldots,N\}$ with $N=\left\lceil\frac1\kappa\right\rceil$. For each $i=1,2$, set $$B_i:=\{a_\ell:\ell \in T_i\}.$$

We claim that the map  $t\mapsto a_t$ is a bijection from $T_i$ onto $B_i$. It suffices to show that this map is injective. In fact, if $a_t=a_s$, then $$|t-s|=|r_s-r_t|<\kappa<\delta_j,$$ and hence the $\delta_j$-separation of $T_i$ implies that the map $t\mapsto a_t$ is injective on $T_i$. The claim is true.

		If $t\in T_1$ and $t'\in T_2$, we write $t=\kappa a_t+r_t$ and $t'=\kappa a_{t'}+r_{t'}$, where $a:=a_t\in B_1$, $a':=a_{t'}\in B_2$ and $0\le r_t,r_{t'}<\kappa$, then
\[tt'=\kappa^2aa'+\kappa a r_{t'}+\kappa a'r_t+r_tr_{t'}.\]
It, together with $\kappa a\le t\le1$, $\kappa a'\le t'\le1$, $|u|\ge L_j\ge1$ and $\varepsilon\le1$, yields that
\[\begin{aligned}
|utt'-\theta aa'|&=\left|u\left(\kappa a r_{t'}+\kappa a'r_t+r_tr_{t'}\right)\right| \le 2|u|\kappa+|u|\kappa^2 			=2\varepsilon+\frac{\varepsilon^2}{|u|} \le 2\varepsilon+\varepsilon^2\le 3\varepsilon.
\end{aligned}\]
Now the elementary inequality $|e^{-2\pi ix}-e^{-2\pi iy}|\le2\pi|x-y|$ gives  \[|e^{-2\pi iutt'}-e^{-2\pi i\theta aa'}|\le6\pi\varepsilon\le\frac{c_0}{2}.\]
By the claim above, we get that $\#B_i=\#T_i=k_i$ for each $i$, and $$\left|\sum_{t\in T_1}\sum_{t'\in T_2} (e^{-2\pi iutt'}-e^{-2\pi i\theta aa'})	\right|\leq  \sum_{t\in T_1}\sum_{t'\in T_2}|e^{-2\pi iutt'}-e^{-2\pi i\theta aa'}|  \le\frac{c_0}{2}k_1k_2,$$
which yields, from \eqref{eq:bilinear-large} and the triangle inequality, that
		\begin{equation*}
			\left|\sum_{a\in B_1}\sum_{a'\in B_2}e^{-2\pi i\theta aa'}\right|\ge \left|\sum_{t\in T_1}\sum_{t'\in T_2}e^{-2\pi iutt'}\right|-\left|\sum_{t\in T_1}\sum_{t'\in T_2}
		(e^{-2\pi iutt'}-e^{-2\pi i\theta aa'})
		\right| \geq \frac{c_0}{2}k_1k_2.
		\end{equation*}
		By Cauchy--Schwarz inequality,
		\begin{align*}
			\left(\frac{c_0}{2}k_1k_2\right)^2 &\leq\left|\sum_{a\in B_1}\sum_{a'\in B_2}e^{-2\pi i\theta aa'}\right|^2 \leq  \left(\sum_{a\in B_1} 1^2 \right)
			\left(\sum_{a\in B_1}  \left|\sum_{a'\in B_2}e^{-2\pi i\theta aa'}\right|^2\right),
		\end{align*} i.e.,
		\begin{equation}\label{eq:theta-size}
			\frac{c_0^2}{4}k_1k_2^2\leq \sum_{a\in B_1}  \left|\sum_{a'\in B_2}e^{-2\pi i\theta aa'}\right|^2.
		\end{equation}
		
We next estimate the right-hand side of \eqref{eq:theta-size} from above.  First, since  $\theta=\frac{\varepsilon^2}{u}$ and $N<\frac{|u|}{\varepsilon}+1,$   it follows that
		\begin{equation*}
			|\theta|N \le\varepsilon+\frac{\varepsilon^2}{|u|} \le\varepsilon+\varepsilon^2 \le\frac14.
		\end{equation*}
		Note that any distinct $a_1,a_2\in B_1$, we have
		\[|\theta|\le|\theta a_1-\theta a_2|
		\le|\theta|N\le\frac14.\]
		Thus, the points $\{x_a:=-\theta a:~a\in B_1\}$ are
		$|\theta|$-separated. Applying the Large sieve (Lemma \ref{lem:large-sieve}) to the degree $N$, separation constant  $\Delta=|\theta|$, and coefficients $c_n=\mathbf 1_{B_2}(n):=\left\{\begin{array}{cc}
                     1, & n\in B_2, \\
                     0, & n\notin B_2,
                   \end{array}
 \right.$ for $0\le n\le N$, we obtain that
\[\sum_{a\in B_1}\left|\sum_{a'\in B_2}e^{-2\pi i\theta aa'}\right|^2\le (N+|\theta|^{-1})k_2.\]
Since $|\theta|^{-1}\leq \frac{N}{\varepsilon}$ and $N\le\frac{|u|}{\varepsilon}+1\le\frac{U_j}{\varepsilon}+1
\le\left(\frac1\varepsilon+1\right)U_j$, we have
 \begin{equation}\label{eq:theta-size1}
\sum_{a\in B_1}\left|\sum_{a'\in B_2}e^{-2\pi i\theta aa'}\right|^2\le \left(\frac1\varepsilon+1\right)^2U_jk_2.
\end{equation}

Thus, by \eqref{eq:theta-size} and \eqref{eq:theta-size1}, we get that
\begin{equation}\label{eq:product-count}
k_1k_2 \le  C U_j,
\end{equation}
where  $C=\frac{4}{c_0^2}\left(\frac1\varepsilon+1\right)^2$  is a constant only depending on $c_0$.  	
		
		Finally, we choose $T_i$ to be a maximal $\delta_j$-separated subset of $E_i$, then it is straightforward to verify that the maximality implies that every point of $E_i$ lies within distance $\delta_j$ of some point of $T_i$. Therefore, \[E_i\subset\bigcup_{t\in T_i}[t-\delta_j,t+\delta_j].\]
As before, we use $N(E,\delta)$ to denote the least number of sets of diameter at most $\delta$ needed to cover $E$. Set \[\tau_j=2\delta_j=\frac{4\varepsilon}{L_j}.\]
Then $N(E_i,\tau_j)\le k_i.$ Consequently, \eqref{eq:product-count} gives \[N(E_1,\tau_j)N(E_2,\tau_j)\le CU_j.\]
		
		Using the inequality $\liminf a_j+\liminf b_j\le\liminf(a_j+b_j)$ and $\tau_j\to0$, the definition of lower box dimension yields
		\[\begin{aligned}
			\underline{\dim}_B E_1+\underline{\dim}_B E_2
			&\le \liminf_{j\to\infty}\frac{\log N(E_1,\tau_j)+\log N(E_2,\tau_j)}{\log(1/\tau_j)} \le\liminf_{j\to\infty}\frac{\log(CU_j)}{\log(L_j/(4\varepsilon))}
			=\liminf_{j\to\infty}\frac{\log U_j}{\log L_j}.
		\end{aligned}\]
		This completes the proof.
	\end{proof}

	\begin{proof}[Proof of Theorem \ref{thm:sharpness}\textup{(b)}]
If $S$ is an $H$-set, then there exist a strictly increasing sequence of positive integers \(q_j\to\infty\) and a nonempty open interval \(I\subset(0,1)\) such that $\{q_jx\}\notin I$ for all $x\in S$ and for all $j\ge1$. By Lemma \ref{lem:avoidance-fourier}, there exist two constants $c_0>0$ and $M\ge1$, only depending on $I$, such that, for every Borel probability measure $\nu$ supported on $S$ and every $j$, there exists an integer $k$ with  $1\leq |k|\leq M$  such that $|\widehat{\nu}(kq_j)|\ge c_0$.
		
		Let $u_j=kq_j$.  Then $q_j\le|u_j|\le Mq_j$. Applying Lemma \ref{lem:fourier-window-product} to $L_j:=q_j$,  $U_j:=Mq_j$, and the two nonempty compact sets $E_1,E_2\subset[0,1]$ satisfying
		$E_1\cdot E_2\subset S$, we obtain that
		\[\underline{\dim}_B E_1+\underline{\dim}_B E_2\le\liminf_{j\to\infty}\frac{\log(Mq_j)}{\log q_j}=1.\]
		Since $\dim_H E\le\underline{\dim}_B E$, this proves the desired assertion.  \end{proof}

The same argument in Theorem \ref{thm.1.8} can be generalized to the general $H^{(n)}$-set.  We now state the standard terminology from \cite[p. 50]{Salem1963}.
\begin{definition}\label{def:hn-set}
Let $n\ge1$. A sequence $\{\mathbf q_j\}_{j\ge1}$,
$\mathbf q_j=(q_j^{(1)},\ldots,q_j^{(n)})\in\mathbb Z^n$, is called \emph{normal} if
\[|\langle \mathbf k, \mathbf q_j\rangle|\to\infty 	\qquad\text{for every }\mathbf k\in\mathbb Z^n\setminus\{0\}.\]
Here, $\langle \cdot,\cdot\rangle$ denotes the inner product of $\mathbb R^n.$	A compact set $S\subset[0,1]$ is called an \emph{$H^{(n)}$-set} if there exist a normal sequence $\{\mathbf q_j\}_{j\ge1}$ and a nonempty open rectangle $V\subset\mathbb T^n$ such that
\[\bigl(\{q_j^{(1)}z\},\ldots,\{q_j^{(n)}z\}\bigr)\notin V
\qquad(z\in S,\ j\ge1).\]
When $n=1$, this reduces to the usual $H$-set condition.
\end{definition}

\begin{proposition}\label{prop:hn-bound}
Let $n\ge1$ and let $S\subset[0,1]$ be an $H^{(n)}$-set associated with a sequence $\mathbf q_j$ and a nonempty open rectangle $V$. If $E_1,E_2\subset[0,1]$ are nonempty compact sets such that $E_1\cdot E_2\subset S$, then there is an integer $M=M(V)\ge1$ such that
\begin{equation}\label{eq:hn-bound}
\underline\dim_B E_1+\underline\dim_B E_2\le 	\liminf_{j\to\infty}\frac{\log U_j}{\log L_j},
\end{equation} where
\[L_j=\min_{\substack{\mathbf k\in\mathbb Z^n\setminus\{0\}\\ \|\mathbf k\|_\infty\le M}}	|\langle \mathbf k, \mathbf q_j\rangle|, \qquad	U_j=	\max_{\substack{\mathbf k\in\mathbb Z^n\setminus\{0\}\\ \|\mathbf k\|_\infty\le M}}	|\langle \mathbf k, \mathbf q_j\rangle|,\] with  $L_j\to\infty$ and $\|\mathbf k\|_\infty=\max\{|k_1|,\ldots,|k_n|\}$.  In particular, if
$$\liminf_{j\to\infty}\frac{\log U_j}{\log L_j}=1,$$
then \[\underline\dim_B E_1+\underline\dim_B E_2\le1.\]
\end{proposition}

\begin{proof}
Let $\psi$ be a nonnegative smooth function  on $\mathbb T^n$, supported in $V$, such that $\int_{\mathbb T^n}\psi=1$.
Expanding \(\psi(\mathbf{x})\) into a multivariate Fourier series
\[\psi(\mathbf{x})=\sum_{\mathbf{k}\in\mathbb{Z}^n}\widehat{\psi}(\mathbf{k})e^{2\pi i\langle\mathbf{k},\mathbf{x}\rangle}\]
where \(\mathbf{k}=(k_1,\ldots,k_n) \in \mathbb{Z}^n\) is the frequency vector, and \(\widehat{\psi}(\mathbf{k})\) are the Fourier coefficients \[\widehat{\psi}(\mathbf{k})=\int_{[0,1]^n}\psi(\mathbf{x})e^{-2\pi i\langle\mathbf{k},\mathbf{x}\rangle}\,d\mathbf{x} \qquad (\mathbf{k}\in\mathbb Z^n).\]
Since the Fourier series of $\psi$ is absolutely convergent, we have $C_\psi=\sum_{\mathbf k\in\mathbb Z^n\setminus\{0\}}|\widehat\psi(\mathbf k)|<\infty.$
Thus, there exists $M\ge1$ so that
\[\sum_{\|\mathbf k\|_\infty>M}	|\widehat\psi(\mathbf k)|\le\frac12.\]

Let $\nu$ be a Borel probability measure supported on $S$. Since $\bigl(\{q_j^{(1)}z\},\ldots,\{q_j^{(n)}z\}\bigr)\notin V
~(z\in S),$	we have
\[0=\int\psi(q_j^{(1)}z,\ldots,q_j^{(n)}z)\,d\nu(z)
=1+\sum_{\mathbf k\ne0}\widehat\psi(\mathbf k)\, 	\widehat\nu(-\langle \mathbf k, \mathbf q_j\rangle).\]
As in the proof of Lemma \ref{lem:avoidance-fourier}, for every $j$ there is $\mathbf k\in\mathbb Z^n\setminus\{0\}$ with $\|\mathbf k\|_\infty\le M$ such that
\[|\widehat\nu(\langle \mathbf k, \mathbf q_j\rangle)|\ge c_0,
\qquad 	c_0=\frac1{2C_\psi}.\]
By the definition of $L_j$ and $U_j$, we have 	$L_j\le |\langle \mathbf k, \mathbf q_j\rangle|\le U_j$ and $L_j\to\infty$. After discarding finitely many indices, we may assume $L_j\ge1$. So the hypotheses of Lemma \ref{lem:fourier-window-product} are satisfied, which therefore gives \[\underline\dim_B E_1+\underline\dim_B E_2\le\liminf_{j\to\infty} \frac{\log U_j}{\log L_j}.\]
The last assertion follows immediately. The proof is complete.
\end{proof}
We conclude this section with a remark on Proposition \ref{prop:hn-bound}. Put $$\rho:=\liminf_{j\to\infty}\frac{\log U_j}{\log L_j}$$ and $E_1=E_2=E$. By  Proposition \ref{prop:hn-bound}, $\underline\dim_B E\le\frac{\rho}{2}$. In particular, if  $L_j$ and $U_j$ have the same logarithmic growth along a subsequence, then $\rho=1$ and hence $\underline\dim_B E\le1/2$. Thus, without further conditions on the frequency scales, the estimate in \eqref{eq:hn-bound} does not necessarily yield the $1/2$ bound for the lower box dimension.  Here we provide two  simple examples.
\begin{example}
If we take \(n=2\), then the sequence \(\mathbf{q}_j=(j,\lfloor\sqrt{2}j\rfloor)\)  normal. Indeed, for every fixed non-zero vector \(\mathbf{k}=(a,b)\in\mathbb{Z}^2\), using the identity \(\lfloor x\rfloor = x - \{x\}\), we have
\[\langle \mathbf{k}, \mathbf{q}_j\rangle = aj + b\lfloor\sqrt{2}j\rfloor = aj + b(\sqrt{2}j - \{\sqrt{2}j\}) = (a+b\sqrt{2})j - b\{\sqrt{2}j\}.
\] Notice that for a fixed bound \(M\), there are only finitely many vectors \(\mathbf{k}\) satisfying \(1\le\|\mathbf{k}\|_\infty\le M\), and   \(|b\{\sqrt{2}j\}| \le |b| \le M\). We obtain that    both \(L_j\) and \(U_j\) are bounded above and below by positive constant multiples of \(j\). Consequently, \(\rho=1\), and Proposition \ref{prop:hn-bound} gives the \(1/2\) lower-box bound.
\end{example}

\begin{example}
If we take \(n=2\), then the sequence  \(\mathbf{q}_j=(j,j^2)\)  is normal. In fact, for any non-zero vector \(\mathbf{k}=(a,b)\in\mathbb{Z}^2\), the inner product is exactly
\[\langle \mathbf{k}, \mathbf{q}_j\rangle = aj + bj^2.\]
Therefore, if  \(1\le\|\mathbf{k}\|_\infty\le M\) for some \(M\ge1\), it is not hard to see that  $U_j=M(j+j^2)$ and $L_j=j$. Consequently, $$\rho = \lim_{j\to\infty} \frac{\log U_j}{\log L_j} = \lim_{j\to\infty} \frac{\log(M(j+j^2))}{\log j} = 2.$$
\end{example}




\end{document}